\documentclass[a4paper,12pt,reqno]{amsart} 
\usepackage{fullpage}
\usepackage{times}
\usepackage{amsfonts}
\usepackage{amssymb}
\usepackage{enumitem}
\usepackage{amsmath}
\usepackage{mathrsfs }
\usepackage{amsthm}
\usepackage{dsfont}
\usepackage{soul}
\usepackage{stmaryrd}

\usepackage{cases}
\usepackage{graphicx} 
\usepackage{calc}
\usepackage[colorlinks=true,linkcolor=black,citecolor=magenta]{hyperref}
\usepackage[dvipsnames]{xcolor}
\usepackage{xifthen} 
\usepackage{appendix}

\newcommand{\VT}[1]{\|#1\|_{\operatorname{TV}}}
\newcommand{\Z}{\mathbb{Z}}
\newcommand{\N}{\mathbb{N}}
\newcommand{\R}{\mathbb{R}}

\newcommand{\E}[1]{\mathbb{E}\left(#1\right)}
\renewcommand{\P}[1]{\mathbb{P}\left(#1\right)}

\newcommand{\Cont}[1]{\operatorname{Cont}\left(#1\right)}
\newcommand{\partt}{\partial_t }
\newcommand{\partx}{\partial_x }

\newcommand{\itg}{\displaystyle\int}
\newcommand{\ind}{\mathds{1}}
\newcommand{\segment}{\llbracket 1,N \rrbracket}
\newcommand{\macroleft}{\ensuremath{\mathscr{L}_p^\alpha}}
\newcommand{\macrostack}{\ensuremath{\mathscr{S}_p^\alpha}}

\newcommand{\ZRP}{\textsc{ZRP }}
\newcommand{\ASEP}{\textsc{ASEP }}

\newcommand{\SEP}{\textsc{SEP }}

\newtheorem{theorem}{Theorem}
\newtheorem{Lemma}{Lemma}
\newtheorem{Proposition}{Proposition}
\newtheorem{corollary}{Corollary}
\newtheorem{definition}{Definition}
\newtheorem{remark}{Remark}
\newtheorem*{assump}{Assumption}

\title{Cutoff Phenomenon for asymmetric Zero Range process with monotone rates }
\author{Ons Rameh }
\address{Ons Rameh: 
DMA, École normale supérieure, Université PSL, CNRS, 75005 Paris, France\\
Université Paris Cité and Sorbonne Université, CNRS, Laboratoire de Probabilités, Statistique et Modélisation, F-75013 Paris, France
}
\email{\href{mailto:orameh@dma.ens.fr}{orameh@dma.ens.fr}}
\date{}
\begin{document}

\begin{abstract}
    We investigate the mixing time of the asymmetric Zero Range process on the segment with a non-decreasing rate. We show that the cutoff holds in the totally asymmetric case with a convex flux, and also with a concave flux if the asymmetry is strong enough. We show that the mixing occurs when the macroscopic system reaches equilibrium.  A key ingredient of the proof, of independent interest, is the hydrodynamic limit for irregular initial data. 
\end{abstract}
\maketitle
\tableofcontents
\section{Introduction}
\subsection{Model}
 The Zero range process on $I_N=\segment$ is a model describing the evolution of $k\geq 1$ indistinguishable particles hopping on the lattice $I_N$ with a rate that depends on the number of the particles at the departure site and the direction of the jump. This model was introduced by Spitzer in 1970 \cite{SPITZER1970246} as a stochastic lattice gas with on-site interaction. Since then, it has been widely studied (see \cite{Evans_2005, 10.1214/aop/1176996977, kipnis1999sli} and the reference therein).\\
Let us now explicitly describe the process. Define
\begin{equation}
    \Omega_{N, k}^0=\left\{\eta \in \N^{\llbracket 1,N \rrbracket}: \sum_{i=1}^N \eta(i)=k\right\}.
\end{equation}
the space of configurations of $k$ particles on $I_N$.\\

Given $p \in (\frac{1}{2},1]$, $q=1-p$ and a rate function by  $g: \N \mapsto [0,+\infty)$ a non-decreasing Lipschitz function such that $g(0)=0< g(1)$,
an Asymmetric Zero Range process (\ZRP$(g,p)$)  is a  continuous time Markov process on $\Omega_{N,k}^0$ with the following generator:
\begin{equation}
    \mathcal{A}_NJ(\eta)= \sum_{i=1} ^ {N-1} p g(\eta(i)) \left(J(\eta^{i,i+1})-J(\eta)\right) +\sum_{i=2} ^ {N} q g(\eta(i)) \left(J(\eta^{i,i-1})-J(\eta)\right) 
\end{equation}
where $\eta^{x,y}$ is the configuration obtained when one particle from $x$ jumps to the site $y$ and $J: \Omega_{N,k}^0 \mapsto \R$  is a test function.

\begin{figure}[h!]
    \label{ZRP: illustration }
    \centering
    \includegraphics[width=10cm]{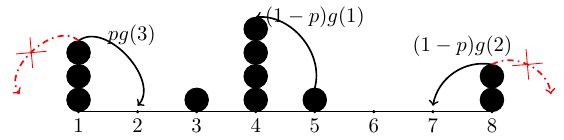}
\caption{A configuration of particles with $k=11$ and $N=8$. The crossed red arrows correspond to forbidden jumps.}
\end{figure}
\begin{remark}
    When $g$ is linear, the dynamic describes the evolution of $k$ independent particles. Each particle jumps to the right with rate $pg(1)$ and to the left with rate $qg(1)$.
\end{remark}
The process is irreducible and reversible with respect to its equilibrium measure 
\begin{equation}
\pi_{N,k}(\eta)=\frac{1}{Z_{N,k}} \prod_{y=1}^N
 \frac{\left(\frac{p}{q}\right)^{\eta(y)y}}{g(\eta(y))!}
\end{equation}
where $g(k)! = g(k)g(k-1)\cdots g(1)$.
In particular, when $p=1$, the equilibrium measure is a Dirac mass on the configuration $\vee_{N,k}= k \delta_N$.

\subsection{Mixing time and Cutoff Phenomenon}Classical Markov process theory guarantees that the process converges to its equilibrium measure $ \pi_{N,k}$ starting from any configuration. This article aims to estimate the speed at which the equilibrium is reached. The standard choice is the total-variation distance, as it allows for uniform control, following Diaconis' earlier work \cite{Diaconis1981-yy, Aldous1986SHUFFLINGCA}. We denote the worst-case total-variation distance by:
\begin{equation}
    d_{N,k}(t)= \max_{\eta \in \Omega_{N,k}^0} \VT{P_t^\eta(\cdot)-\pi_{N,k}}
\end{equation}
where $P_t^\eta$ is the law of a \ZRP at time $t$ started at $\eta$.

For any level of precision $\theta>0$, we define the  mixing time of \ZRP: 
\begin{equation}
    T_{\operatorname{mix}}^ {N,k}(\theta):= \inf\{ t\geq 0; \quad d_{N,k}(t)\leq \theta\}.
\end{equation}
In particular, we aim to describe the behaviour of the mixing time when the size of the system and the number of particles diverge to infinity. In this article, we work in the natural setting for a scaling limit, where $\frac{k}{N}$ converge to $\alpha \in \R_+$, a macroscopic finite density. \\

Of particular interest, is the dependence of the mixing time on the precision $\varepsilon$. We say that the system exhibits \emph{a Cutoff phenomenon} if 
$$
\forall \theta \in (0,1), \qquad \frac{T_{\operatorname{mix}}^ {N,k}(\theta)}{T_{\operatorname{mix}}^ {N,k}(1/4)} \xrightarrow[N \rightarrow +\infty]{} 1.
$$
In other words, the distance to the equilibrium stays around $1$ and drops abruptly to $0$. This remarkable phase transition was discovered in the context of card shuffling \cite{Diaconis1981-yy, Aldous1986SHUFFLINGCA}.  We refer the reader to \cite[Chapter $18$]{LevinPeresWilmer2006} for more background on the Cutoff Phenomenon.
\subsection{Previous Results}
    The Cutoff was investigated in various contexts, among which the simple exclusion process that describes the evolution of particles hopping on a lattice and constrained to not coexist on the same site.  \textsc{SEP}'s mixing time has been studied in both symmetric  \cite{Wilson_2004, Lacoin2013MixingTA,Lacoin_2016, lacoin2016simpleexclusionprocesscircle} and  (weakly) asymmetric \cite{Levin_2016,benjamini2002mixingtimesbiasedcard, Labb__2019, labbé2018mixingtimecutoffweakly} settings. The articles \cite{Lacoin_2016,Labb__2019,labbé2018mixingtimecutoffweakly,labbé2022hydrodynamiclimitcutoffbiased} underlined a notable link between the macroscopic relaxation to equilibrium (the Hydrodynamic limit) and the mixing time. 
    In particular, it was proven that in the asymmetric setup, the macroscopic equilibrium is reached in a finite time and it corresponds exactly to the mixing time.  This link was also emphasised for the symmetric \SEP with reservoir \cite{gonçalves2021sharpconvergenceequilibriumssep} where the study of hydrodynamic limit allows for a precise description of the Cutoff profile.
    
\ZRP can be seen as a generalization of \textsc{SEP} as it can be mapped into a simple exclusion process (\textsc{SEP}) when the rate function $g$ is constant.  Cutoff for \ZRP was also shown when the rate function is constant on the complete graph \cite{merle2018cutoffmeanfieldzerorangeprocess}. For more general rate functions, it was only studied to the best of our knowledge on the complete graph \cite{hermon2018cutoffmeanfieldzerorangeprocess, tran2022meanfieldzerorangeprocessunbounded}. 

In the context of asymmetric \textsc{ZRP}, the hydrodynamic limit was thoroughly studied on $\Z^d$ \cite{Rezakhanlou1991HydrodynamicLF} for regular initial data. It was also established in the non-conservative case \cite{bahadoran2011hydrodynamics} and in a random environment \cite{bahadoran2020zerorangeprocessrandomenvironment}.  {In all these results the rate function is assumed nondecreasing to guarantee the attractivity of the process and to avoid the condensation phenomena (see for example \cite{condensation}). }

\subsection{Our contribution} 
The main achievement of this paper is to prove that cutoff holds in the totally asymmetric convex flux, and in the concave flux if the asymmetry is strong enough (Theorems \ref{thm lower bound mixing time } and \ref{thm upper bounds mixing time }). Moreover, we show that it occurs exactly at the macroscopic equilibrium time, which highlights once again the link between hydrodynamic limit and mixing times. Our proof exploits the hydrodynamic limit of \ZRP on the segment, which we generalise for irregular initial data (Theorem \ref{thm hydro limit on the segment}). As far as we know, the hydrodynamic limit for irregular initial data hasn't been studied before. We also give a precise description of the evolution of the position of the left-most particles and the stack of particles at position $N$.
\subsection{Acknowledgment} The author warmly thanks Max Fathi for his advice throughout this project and for thoroughly reading the drafts. The author also kindly thanks Christophe Bahadoran, Cyril Imbert, Cyril Labbé, Hubert Lacoin, Fraydoun Rezakhanlou, Ellen Saada, Justin Salez, Marielle Simon and Cristina Toninelli for helpful discussions. This work received funding from the Agence Nationale de la Recherche (ANR) Grant
ANR-23-CE40-0003 (Project CONVIVIALITY), as well as funding from the Institut Universitaire de France.
\section{Main Results}
Before we state the main theorems, we introduce some notation. For any $\varphi < \varphi ^\star$, we define:
$$
Z(\varphi)=\sum_{k \geq 0} \frac{\varphi^k}{g(k)!} \qquad \text{and} \qquad \forall k \in \N, \, \bar{\nu}_{\varphi}(k)=\frac{1}{Z(\varphi)} \frac{\varphi^k}{g(k)!}
$$
where $\varphi^\star$ is the radius of convergence of $Z$. Throughout the paper, we will always assume the following:
\begin{assump} \hfill
    \begin{itemize}
        \item  $g$ is non--decreasing,
        \item $g$ is Lipshitz :
       $$
g^\star :=\sup _{k \geq 0}|g(k+1)-g(k)|< \infty,
$$
        \item $
    \lim _{\varphi \uparrow \varphi^*} Z(\varphi)=\infty.$
    \end{itemize}
\end{assump}
The first assumption ensures the attractivity of the process and avoids potential condensation phenomena that may arise. 
The family of measure $(\bar \nu_ \varphi)_\varphi$ can be thus parametrized by its expectation: 
$$\forall \alpha \in \R_ +, \, \nu (\alpha) = \bar \nu_{\Phi(\alpha)} $$
where $\Phi$ is the inverse function of $\varphi \rightarrow \mathbb{E}_{\bar{\nu}_{\varphi}}\left(X\right)$ and corresponds to the macroscopic flux for a density $\alpha$.
When $\Phi$ is convex, we also denote by $\Psi$ its convex conjugate:
$$\Psi(x) = \sup_{\alpha}(\alpha x-\Phi(\alpha)).$$
If $\Phi$ is concave, $\Psi$ denotes the convex conjugate of $(-\Phi)$ instead. We recall that the convexity (resp. concavity) of $g$ implies that $\Phi$ is strictly convex (resp. strictly concave) \cite[Proposition 3.1]{balazs2007convexitypropertyexpectationsexponential}.\\

Let $\mathcal{M}_{+}([0,1])$ be the space of non--negative measure on $[0,1]$ endowed with the weak topology. For any  sequence of configurations
 $(\eta_{N,0})_N \in \Omega_{N,k}^0$ ,  we define the empirical measure associated with the process $(\eta_{N,t})$:
$$ \rho_{N,t} (dy) = \frac{1}{N} \sum_{i \in I_N}  \eta_{N,t}(i) \delta_{i/N}(dy) \in \mathcal{M}_{+}([0,1])$$
where $(\eta_{N,t})$ is \ZRP started at $(\eta_{N,0})_N$. Notice that $(\rho_{N,t})_t$ belongs to the  $\mathbb{D}(\R_+, \mathcal{M}_{+}([0,1]))$ the set of cadlag functions with values in $\mathcal{M}_+([0,1])$ endowed with the Skorokhod topology. This allows us to simultaneously embed the configurations with different sizes in a common space.   \\

In the context of \ZRP, it is intuitively clear that the worst initial condition should correspond to $\wedge_{N,k}=k \delta_1$ where all particles are on the first site. From a macroscopic scale, it corresponds to imposing $ \alpha \delta_0$ as initial data where $\alpha$ would be the macroscopically observed mass.  The following theorem gives the hydrodynamic limit of \ZRP on the segment for general irregular initial data.
\begin{theorem}
    \label{thm hydro limit on the segment}
    Assume that $\Phi$ is either strictly convex or strictly concave.
    Given  $u_0 \in \mathcal{M}_{+}([0,1]) $ with total mass $\alpha$, let $U$ be the unique Barron--Jensen viscosity solution to
    \begin{equation}
 \tag{HJ}
\label{eq: HJ on R}
\left\{\begin{array}{l}
\partt U+ \Phi(\partx U)=0, \quad t \geq 0, x \in \R, \\
U(0,\cdot)= \itg_ {[0,\cdot)} u_ 0(dy).
\end{array}\right.
\end{equation}
For any $p>\frac{1}{2}$, if $ (\rho_N(0, dy))_N $  weakly converge to $u_0(dx)$,
    the sequence $(t\rightarrow \rho_{N,\frac{Nt}{p-q}})_N$ converge in law in $\mathbb{D}(\R_+, \mathcal{M}_+(\R))$ to a measure-valued function $t \rightarrow v(t,dx)$.
    Moreover, $$\forall t>0, \forall x \in (0,1], \, v(t,[0,x])=\max\left(U, \alpha \ind_{\R_+ \times[1, +\infty)}\right)^\star. $$
\end{theorem}

Let us underline that we chose to speed up the process by $\frac{N}{p-q}$ to ensure the limit does not depend on $p$. The precise definition of Barron--Jensen viscosity solution of the Hamilton--Jacobi equation (\ref{eq: HJ on R}) is delayed to subsection \ref{subseq: HJ def}. The limit is usually determined (see \cite{Rezakhanlou1995MicroscopicSO,Evans_2005}) as the unique entropic solution to  the scalar conservation law:
\begin{equation}
\label{eq: scalar csv}
\tag{CL}
\left\{\begin{array}{l}
\partt u+ \partx \Phi(u)=0, \quad t\geq 0, x \in \mathbb{R}, \\
u(0,\cdot)=u_0.
\end{array}\right.
\end{equation}
Both notions of solutions are equivalent when $u_0 \in L^1 \cap L^\infty$ for convex or concave flux $\Phi$ \cite{csv_law_equivalence}. Beyond the difficulties raised by the presence of boundaries, the irregularity of initial data with eventually non-bounded rate jumps required the use of Hamilton--Jacobi equations to uniquely identify the limit $v$. Note that
although there exists a notion of measure-valued solution of scalar conservation law, ensuring well-posedness in a more general framework (see \cite{demengelserredoi:10.1080/03605309108820758,bertsch2018uniquenesscriterionmeasurevaluedsolutions,Bertsch2018RadonMS,Bertsch2020DiscontinuousSO}), adapting the uniqueness criterion to our situation seems more challenging. Regarding the boundary conditions, we extensively use the attractivity of the process and the fact that the characteristics point towards the right. 
\\

For any $\alpha \geq 0$, let $U^\alpha$ be the Barron--Jensen viscosity solution associated to $u_{0}^\alpha= \alpha \delta_ 0 $.
One can in particular remark that the macroscopic system reaches equilibrium in finite time. We define for any $\alpha \geq 0$ the macroscopic equilibrium time
\begin{equation}
    T_{eq,p,\alpha}= \frac{1}{p-q} \inf\{t >0, U^\alpha(t,dx) \leq \alpha \ind_{x>1}\}.
\end{equation}
If the macroscopic system has not reached the equilibrium, the microscopic system is also still far from it. Consequently, the macroscopic equilibrium time provides the following lower--bound on the mixing time. 

\begin{theorem}
    \label{thm lower bound mixing time }
    Assume that $\Phi$ is either strictly convex or strictly concave. For any $p > \frac{1}{2}$,
    $$
    \forall \alpha \in \R_+, \; \forall \theta \in (0,1),\qquad \liminf _{\substack{N \rightarrow \infty \\ k / N \rightarrow \alpha}} \frac{T_{\text {mix }}^{N, k}(\theta)}{N} \geq T_{eq,p,\alpha}.
    $$
\end{theorem}
Our main result is that the lower--bound is sharp in the totally asymmetric convex case and the concave setting if the asymmetry is strong enough.
\begin{theorem}\hfill
    \label{thm upper bounds mixing time }
    \begin{enumerate}[ref=\ref{thm upper bounds mixing time }.\arabic*]
        \item \label{thm upper bdd cvx} For $p=1$ and assuming that $\Phi$ is strictly convex and $g^\star=g(1)$,
         $$
    \forall \alpha \in \R_+, \, \forall \theta \in (0,1), \quad \lim _{\substack{N \rightarrow \infty \\ k / N \rightarrow \alpha}} \frac{T_{\text {mix }}^{N, k}(\theta)}{N} =T_{eq,1,\alpha}= \frac{1}{g(1)}.
    $$
        \item \label{thm upper bdd concave} Assume that $\Phi$ is strictly concave and bounded by $\bar g$.  Assume $p>\frac{1}{2},\alpha>0$  are such that
        \begin{equation}
            \label{cond sur p et alpha}
            \frac{pg(1)-q  \bar g}{p-q} > \max\left\{\frac{\Phi\left(\Psi^\prime \left(\Psi^{-1}\left(\frac{\alpha}{T_{eq,1,\alpha}}\right)\right)\right)}{\Psi^\prime \left(\Psi^{-1}\left(\frac{\alpha}{T_{eq,1,\alpha}}\right)\right)} ,\Phi\left((-\Phi^\prime)^{-1}\left(\frac{-1}{T_{eq,1,\alpha}}\right)\right) \right \}.
        \end{equation}
        Then,
        $$\forall \alpha \in \R_+, \, \forall \theta \in (0,1), \quad \lim _{\substack{N \rightarrow \infty \\ k / N \rightarrow \alpha}} \frac{T_{\text {mix }}^{N, k}(\theta)}{N} =T_{\operatorname{eq},p, \alpha}.$$
    \end{enumerate}
\end{theorem}
Clearly, it implies that the \ZRP in this framework exhibits a Cutoff since the limit does not depend on $\theta$. As was already the case in the \ASEP \cite{Labb__2019}, the hydrodynamic limit does not provide precise enough information to prove the matching upper--bound. An easy way to see it is by considering $1 \delta_1 + (k-1) \delta_N$ as an initial condition where macroscopic equilibrium is achieved and microscopically we need to wait for $(p-q)g(1)N$.\\
The main idea is to control precisely:
\begin{itemize}
    \item $L_{N,k}$: the position of the left-most particle 
    \item $S_{N,k}$: the number of particles at the position $N$
\end{itemize}
and prove that they do follow their macroscopic analogue $\left(\macroleft,\macrostack\right) $, predicted by the hydrodynamic limit. Once these two statistics reach their equilibrium state, the system rapidly mixes. At this level of generality, $\pi_{N,k}$, $\macroleft$ and $\macrostack$ are not explicit, making the problem delicate. We also emphasise that \ZRP do not have a particle-hole duality, unlike \ASEP. Thus, the study of $L_{N,k}$ and $S_{N,k}$ need to be conducted separately. The technical condition (\ref{cond sur p et alpha}) is due to a comparison to \ASEP case which lacks precision. Indeed, condition (\ref{cond sur p et alpha}) requires that the maximum macroscopic speed and flux achieved by the \ZRP be below the critical level $pg(1)-q\bar g$ of the \ASEP, to which we are comparing (see remarks \ref{rem: cond  explication 1} and \ref{rem cond explication 2}).  { We believe that the result should hold in a more general context.} 
\begin{remark}
The condition (\ref{cond sur p et alpha}) is verified on a wider class of \ZRP than just for the constant rate. Indeed, we observe that when $\alpha$ converge to $0$,
  $\Phi\left((-\Phi^\prime)^{-1}\left(\frac{-1}{T_{eq,1,\alpha}}\right)\right) \rightarrow 0$. Hence, one can fix $\alpha_{\max}$ such that 
  $$\Phi\left((-\Phi^\prime)^{-1}\left(\frac{-1}{T_{eq,1,\alpha}}\right)\right) <g(1).$$
  Thus, when $p=1$, the condition (\ref{cond sur p et alpha}) is trivially fulfilled for any $\alpha \leq \alpha_{\max}$. Consequently, $$ \exists \alpha_{\max}, \; \forall 0< \alpha \leq \alpha_{\max}, \; \exists p_{\min} \in(1/2,1),\; \forall p  \geq p_{\min} \quad  \text{such that the condition (\ref{cond sur p et alpha}) is fulfilled.}$$ 
\end{remark}
\subsection{Organisation of the paper} The rest of the article is organized as follows. In Subsections \ref{subsec: couplage} and \ref{subsec: link with exclusion}, we recall the attractivity property of \ZRP and the mapping to the exclusion process when the rate is constant. Then, we recall the definition of the Barron--Jensen viscosity solution in Subsection \ref{subseq: HJ def} and generalize the hydrodynamic limit on $\Z$ with irregular initial data in Subsection \ref{subsec: hydro limit on Z}. In Section \ref{sec: hydro limit on the segment}, we prove the hydrodynamic limit of \ZRP on the segment. The proof of Theorems \ref{thm lower bound mixing time } and \ref{thm upper bounds mixing time } occupies the whole Section \ref{sec: mixing time proof}. Section \ref{sec: control de  L et de S} is devoted to the control of the position of the left-most particle and the stack at position $N$ when the flux is concave. A list of symbols used is inducted at the end of the Appendix.
\section{Preliminaries}
\label{sec:  Preliminaries}
\subsection{Coupling, attractivity, and orders}
\label{subsec: couplage}
Throughout the article, we shall frequently need to couple \ZRP on different graphs and rate functions. For convenience, we introduce a generic form that serves as a standard framework, to which all specific cases can be reduced.  Given $g_+,g_-: \N \times \Z \mapsto \R_ +$ , we define the following generator: 
\begin{equation}
    \label{eq: generic form generator}
    \mathcal{A}_{g_+,g_ -}J(\eta)= \sum_{i \in Z} g_+(\eta(i),i) \left(J(\eta^{i,i+1})-J(\eta)\right)  +  g_-(\eta(i),i) \left(J(\eta^{i,i-1})-J(\eta)\right)  ,
\end{equation}
where $g_+$ (resp. $g_-$) is the rate of jumps to the right (resp. left) from position $i$. Let us mention that any \ZRP generator on any sub-graph of $\Z$ could be written in this generic form.\\

It can be more convenient to compare their evolution by representing a configuration $\eta$ through its height function. For any $\eta \in \N ^\Z$, we define the height function by 
\begin{equation}
    h_\eta(j)= \sum_{i=1}^{j} \eta(i), \qquad \forall j \in \Z.
\end{equation}
We denote by $\Omega_\Z$ the set of height functions corresponding to configurations on $\Z$.
It is endowed  with the order defined by
$$
\left\{h \geq \Tilde{h}\right\} \Leftrightarrow\left\{\forall i \in \Z, h(i) \geq \Tilde{h}(i)\right\}.
$$
 As the asymmetry induces a bias to the right, we can understand this partial order as a comparison of the evolution of the processes, where ${h}$ is farther than $\tilde h$ from the equilibrium. 
\begin{remark}
    In the case of \ZRP on the segment, we will denote the set of height functions by $\Omega_{N,k}= h(\Omega_{N,k}^0)$. In particular, the configuration $\vee_{N,k}=k \delta_N$ is the minimal element and $\wedge_{N,k}= k \delta_1$ to the maximal element. Due to the bias to the right, the configuration $\wedge_{N,k}$ corresponds to the worst initial condition.
\end{remark}

The following lemma provides a straightforward criterion to compare different processes, which is a slight generalization of the usual attractive coupling. 

\begin{Lemma}[Attractivity]
    \label{lem:attractivity}
    Given $g_+,g_-, \tilde g_+,\tilde g_-: \N \times \Z \mapsto \R_+$ such that:
    \begin{equation}
        \label{cond: attractivity}
        \forall i \in \Z, \; \forall l \geq m, \qquad \tilde g_ +(l,i) \geq g_+(m,i) \text{ and } \tilde g_ -(m,i) \leq g_-(l,i).
    \end{equation}
    For any $\tilde h_0 \leq  h_0 \in \Omega_\Z$, there exists a coupling such that almost surely:
   $$\forall t\geq 0 , \; \tilde h_t \leq  h _t$$
   where $h$ (resp. $\Tilde{h}$) follows the dynamic $\mathcal{A}_{g_+,g_ -}$ (resp. $\mathcal{A}_{\tilde g_+, \tilde g_ -}$).
\end{Lemma}
\begin{proof}
    We introduce a coupling via the following Markov generator:
    \begin{equation}
        \label{eq: coupling attractivity}
        \tag{$\mathscr{C}$}
            \Tilde{\mathcal{A}} J({\eta}, \tilde{\eta}) =
            \sum_{i \in \Z} \mathcal{G}_{i, +} + \mathcal{G}_{i, -}
    \end{equation}
    where $$
\begin{aligned}
\mathcal{G}_{k, \pm}J({\eta}, \tilde{\eta}) &=\left(g_\pm(\eta(i),i) \wedge \tilde g_\pm(\tilde \eta(i),i)\right)\left(J\left(\eta^{i, i \pm 1 }, \tilde \eta^{i,i\pm 1 }\right)-J(\eta, \tilde \eta)\right) \\
& +\left[g_\pm(\eta(i),i)-\tilde g_\pm(\tilde \eta(i),i)\right]_ +\left(J\left(\eta^{i, i \pm 1 }, \tilde \eta\right)-J(\eta, \tilde \eta)\right) \\
& +\left[\tilde g_\pm(\tilde \eta(i),i)-g_\pm(\eta(i),i) \right]_ + \left(J\left(\eta , \tilde \eta^{i,i\pm 1 }\right)-J(\eta, \tilde \eta)\right). 
\end{aligned}$$
As usual, we want particles for the two processes to jump simultaneously as much as possible.  
We denote by $$\mathcal{E}=\{(h,\Tilde{h}) \in \Omega_\Z^2, \;  h\geq \tilde h \}.$$ It remains to show that if $h \geq \tilde h $ then $ \Tilde {\mathcal{A}}\ind_{\mathcal E}(h,\tilde h)\geq 0 $. 
We remark that the jump from $i \in \Z$ that may contribute negatively are:
\begin{itemize}
    \item A jump to the right of an  $\eta$--particle without an $\tilde \eta$--particle jump when $h(i)=\tilde h(i)$. In such a case, $\eta(i) \leq \tilde \eta (i)$ which implies that $g_ +(\eta(i),i)\leq \tilde g_ +(\tilde \eta (i),i)$. 
    \item A jump to the left of an  $\tilde \eta$--particle without an $ \eta$--particle jump when $h(i-1)=\tilde h(i-1)$. In such a case, $\eta(i) \geq \tilde \eta (i)$ which implies that $g_ -(\eta(i),i)\geq \tilde g_ -(\tilde \eta (i),i)$. 
\end{itemize}
In both cases, the corresponding term in the generator is null which concludes the proof.
\end{proof}
\begin{remark}
    As $g$ is non--decreasing, we recover the attractivity of \ZRP :\\
    \indent For any $\tilde h_0\leq h_0$, there exists a coupling such that almost surely for any $t \geq 0$,  $\tilde  h_t\leq h_t$.\\
    One can see the attractivity as the analogue of the $L^1$--contraction of entropy solution to (\ref{eq: scalar csv}). It can also be seen as a microscopic counterpart of the Hamilton--Jacobi's maximum principle. The attractivity is the key argument in the proof of the hydrodynamic limit of the \ZRP in \cite{Rezakhanlou1991HydrodynamicLF},\cite[Chapter 6]{kipnis1999sli} and also for Theorem \ref{thm hydro limit on the segment}.
\end{remark}

We will also seek to transfer information about the macroscopic behaviour of \ZRP on $\mathbb{Z}$ to the ZRP on subintervals. In the spirit of \cite[Lemma $18$]{Labb__2019} and \cite[Lemma 3.3]{bahadoran2011hydrodynamics}, the same coupling (\ref{eq: coupling attractivity}) guarantees a comparison of both processes for short times and far from the border. The proof is identical that in the references but we include it for completeness.
\begin{Lemma}[Comparison far from borders]
    \label{lem local control in coupling ray - line 1}
     Given $a< b \in \Z\cup \{\pm \infty\}, \beta \in\R_ +$, let $\eta$ be a \ZRP$(g,p)$ on $\Z$ and $\Tilde{\eta}$ a \ZRP$(g,p)$ on $\llbracket aN,bN \rrbracket $ coupled according to (\ref{eq: coupling attractivity}).\\
     There exists a constant $C>0$ such that for any $ a^\prime , b^\prime \in \R$ such that $a<  a ^\prime < b ^\prime < b$
       $$
\lim _{N \rightarrow \infty} \sup _{\substack{{\eta}_0=\tilde{\eta}_0 \in \N^{\Z }\\ \sum_i \Tilde{\eta}_0(i) \leq \beta N}}\sup _{s \in[0, C \varepsilon]} \mathbb{E}\left[\frac{1}{N} \sum_{i = a^\prime N }^{b^\prime N}\left| \eta_{\frac{sN}{p-q}}( i)-\tilde{\eta}_{\frac{sN}{p-q}}( i)\right |\right]=0,
$$
where $\varepsilon= \min(b-b^\prime, a^\prime-a)$.
\end{Lemma}

\begin{proof}[Proof of Lemma \ref{lem local control in coupling ray - line 1}]
Given such $ a^\prime , b^\prime \in \R$. We introduce
\begin{itemize}
    \item $H:\mathbb{R} \mapsto [0,1]$ smooth, non-increasing, equal to 1 on $\mathbb{R}_{-}$and to 0 on $[1, \infty)$.
    \item $q: \mathbb{R} \mapsto[0,1]$ equal to 0 on $\left[a^ \prime , b^ \prime \right ]$, to 1 on $(-\infty, a^ \prime - \frac{2\varepsilon}{3}] \cup[b^ \prime +\frac{2\varepsilon}{3}, \infty)$ and such that $\left\|q^{\prime}\right\|_{\infty} \leq 3 \varepsilon^{-1}$.
    \end{itemize} 
We denote by $J(t,x)=H\left(\frac{6g^\star}{(p-q)\varepsilon}t+q(x)\right)$. Then, there exists $c,C>0$  such that $J(t,x)> c$ for any $(t, x) \in[0, C \varepsilon] \times[a^\prime, b^\prime]$ which implies: 
\begin{align*}
    \E{\frac{1}{N} \sum_{i = a^\prime N }^{b^\prime N}\left| \eta_{\frac{tN}{p-q}}( i)-\tilde{\eta}_{\frac{tN}{p-q}}( i)\right |} \leq \E{\frac{1}{cN} \sum_{i \in \Z} J\left(t,\frac{i}{N}\right) \left|\eta_{\frac{tN}{p-q}}(i)-\tilde\eta_{\frac{tN}{p-q}}(i)\right|}.
\end{align*} 
Since ${\eta}_0=\tilde{\eta}_0 $, we have:
\begin{equation}
    \begin{aligned}
        \sum_{i \in \Z} J\left(t,\frac{i}{N}\right) \left|\eta_{\frac{tN}{p-q}}(i)-\tilde\eta_{\frac{tN}{p-q}}(i)\right|= M_t + I_t
    \end{aligned}
\end{equation}
where $M$ is a martingale with $M_0=0$ and 
$$
I_t = \itg_0^{t}ds \left[ \sum_{i \in \Z} \partt J\left(s,\frac{i}{N}\right)  \left|\eta_{\frac{sN}{p-q}}(i)-\tilde\eta_{\frac{sN}{p-q}}(i)\right| + \frac{N}{p-q}J\left(s,\frac{i}{N}\right) \tilde{\mathcal{A}} \left|\eta_{\frac{sN}{p-q}}(i)-\tilde\eta_{\frac{sN}{p-q}}(i)\right|\right].$$
As $J$ is equal to $0$ as soon as $x \notin [a^ \prime - \frac{2\varepsilon}{3}, b^ \prime +\frac{2\varepsilon}{3}]$, we obtain:
\begin{align*}
     \sum_{i \in \Z} & J\left(\frac{s}{N},\frac{i}{N}\right) \tilde{\mathcal{A}} \left|\eta_s(i)-\tilde\eta_s(i)\right| \\
     &\leq \sum_{i \in \Z} p \left[J\left(s,\frac{i+1}{N}\right)-J\left(s,\frac{i}{N}\right)\right] \left|\eta_{\frac{sN}{p-q}}(i)-\tilde\eta_{\frac{sN}{p-q}}(i)\right| \operatorname{ord}_{i,i+1}\\
      & \qquad+  \sum_{i \in \Z} q \left[J\left(s,\frac{i-1}{N}\right)-J\left(s,\frac{i}{N}\right)\right] \left|\eta_{\frac{sN}{p-q}}(i)-\tilde\eta_{\frac{sN}{p-q}}(i)\right|  \operatorname{ord}_{i,i-1},
\end{align*}
where $\operatorname{ord}_{i,j}= \ind_{(\eta(i)-\tilde \eta(i))(\eta(j)-\tilde \eta(j)) \geq 0}$. Since the number of particles is bounded by $\beta N$ and $J$ is smooth, there exists $C_{J,\beta}>0$ such that for any $t\leq C \varepsilon$
\begin{equation*}
    \begin{aligned}
        I_t \leq C_{J,\beta}+ \itg_0^{t}ds \sum_{i \in \Z} (\partt J+ \frac{2 g ^\star}{p-q} |\partx J|)   \left(\frac{s}{N},\frac{i}{N}\right)  \left|\eta_{\frac{sN}{p-q}}(i)-\tilde\eta_{\frac{sN}{p-q}}(i)\right| 
    \end{aligned}.
\end{equation*}
Finally, we observe that
$
\partt J  + \frac{2 g ^\star}{p-q} |\partx J| \leq 0
$
which implies for any $t\in [0, C \varepsilon]$: 
\begin{align*}
    \E{\frac{1}{cN} \sum_{i \in \Z} J\left(t,\frac{i}{N}\right) \left|\eta_{\frac{tN}{p-q}}(i)-\tilde\eta_{\frac{tN}{p-q}}(i)\right|}\leq \frac{C_{J,\beta}}{N},
\end{align*}
and concludes the proof.
\end{proof}

\subsection{Link with the exclusion process}
\label{subsec: link with exclusion}
      Let us recall the generator of the asymmetric exclusion process  (\ASEP) on a segment: 
    \begin{equation*}
    \mathcal{A}_N^{\operatorname{ex}}J(\eta)= \sum_{i=1} ^ {N-1} p  \eta(i)[1-\eta(i+1)]  \left(J(\eta^{i,i+1})-J(\eta)\right) +\sum_{i=2} ^ {N} q  \eta(i)[1-\eta(i-1)]  \left(J(\eta^{i,i-1})-J(\eta)\right) .
\end{equation*}
    For any configuration $\xi \in \Omega_{N,k}^{0,\operatorname{ex}}=\left\{ \xi\in \{0,1\}^{I_{N+k-1}}, \; \sum_{i=1}^{N+k-1}\xi_i= k\right\}$, we define $E(\xi) \in \Omega_{N,k}^0$ according the following rules:
    \begin{itemize}
        \item Label $(z_i)_{i=1}^N$ the position of empty-sites on the configuration $\xi$. Fix $z_0=0$ and $z_N=N+k+1$
        \item On each site $i$,  put $z_i-z_{i-1}-1$ particles.
    \end{itemize}
    \begin{figure}[h!]
    \label{ZRP to exclusion: illustration }
    \centering
    \includegraphics{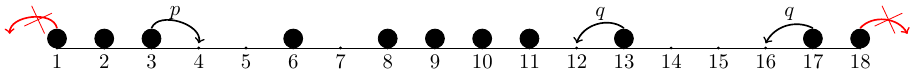}
\caption{The exclusion configuration corresponding to the figure (Fig.\ref{ZRP: illustration })}
\end{figure}
It is easy to see that $E$ is a bijection from $\Omega_{N,k}^{0,\operatorname{ex}}$  to $\Omega_{N,k}^{0}$. Moreover, it is well-known that \ZRP is in bijection with \ASEP on a segment $I_{N+k-1}$ with $k$ particles if the rate function is constant \cite{Evans_2005,Dandekar_2022}. 
\begin{Lemma}[Link between \ZRP and \textsc{SEP}]
    \label{Lemma: link ZRP-exclusion}
     If $(\eta_t)_{t\geq 0}$ is a \ZRP with a rate function $g=\ind_{x\geq 1}$ started at $\eta_0$, then process $\left(E^{-1}(\eta_t)\right)_t$ evolves according to  $\mathcal{A}_N^{\operatorname{ex}}$ started at $E^{-1}(\eta_0)$.
\end{Lemma}

\subsection{Hamilton Jacobi Equation}
\label{subseq: HJ def}
In this subsection, we recall uniqueness and stability properties of the Hamilton-Jacobi equation  on $\Bar{\Omega} = [0, +\infty) \times \R$ 
 \begin{equation}
 \tag{\ref{eq: HJ on R}}
\left\{\begin{array}{l}
\partt U+ \Phi(\partx U)=0, \quad t,x \in \Omega \\
U(0,\cdot)=U_0,
\end{array}\right.
\end{equation}
where $\Omega = (0, +\infty) \times \R$.  Let us assume that $\Phi$ is strictly convex. 

Barron and Jensen introduced in  \cite{BarronJensen} a notion of viscosity solution adapted for discontinuous initial data. 
\begin{definition}
A bounded lower semicontinuous function is a Barron--Jensen solution (\ref{eq: HJ solution discontinuous}) iff for all $\varphi \in \mathcal{C}^1(\Bar{\Omega})$ and at each maximum point $(t_0,x_0) \in \Omega $ of  $U- \varphi$, we have 
         \begin{equation}
                \label{eq: HJ solution discontinuous}
                \tag{BJ}
              \partt \varphi(t_0,x_0) + \Phi(\partx \varphi (t_0,x_0)) =0 .
         \end{equation}
\end{definition}

\begin{remark}
    The Barron--Jensen notion of viscosity solution is also called the ``bilateral solution'' and is equivalent to the Crandall--Lions viscosity solution if $\Phi$ is convex and   $U \in W_{\text {loc }}^{1, \infty}(\Omega)$ \cite[Theorem 2.3]{BarronJensen}.
\end{remark}
This notion ensures the uniqueness of the solution to (\ref{eq: HJ on R}).
\begin{theorem}[{\cite[Theorem 5.14]{Barles1994SolutionsDV}}]
    \label{Theorem uniqueness HJ}
    Assume that $U_0$ is a l.s.c bounded function. There exists at most one l.s.c function $U$ which satisfies both (\ref{eq: HJ solution discontinuous}) and
$$
u(x, 0)=u_0 \quad \text { on } \mathbb{R}^N.
$$
\end{theorem}
It is also well-known that if the initial data $U_0$ is Lipshitz, the viscosity solution is given by the Hopf--Lax formula:
\begin{equation}
    \label{eq: hopf lax formula cvx}
    U(t,x)=\inf _{y \in \mathbb{R}}\left(U_0(y)+t \Psi\left(\frac{x-y}{t}\right)\right)
\end{equation}
where we recall that $\Psi$ is the convex conjugate function of $\Phi$.
Demengel and Serre proved in \cite{demengelserredoi:10.1080/03605309108820758} that the Hopf--Lax formula is stable.
\begin{Proposition}[{\cite[Proposition 4.1]{demengelserredoi:10.1080/03605309108820758}}]
    \label{prop cvg formule de LAX}
    Given a sequence $(u_{0,n})$  bounded in $L^1(\R)$ and converging vaguely to a measure $u_{0,\infty}$, let
    $$\forall n\in \N \cup \{\infty\}, \qquad U_n(t,x)=\inf _{y \in \mathbb{R}}\left(\itg_{(-\infty,y)} u_{0,n}(z)dz+t \Psi\left(\frac{x-y}{t}\right)\right).
$$
Then, for any $t\geq 0, x \in \operatorname{Cont}(U(t))$, $(U_n(t,x))$ converge to $U_\infty(t,x)$ defined by the Hopf--Lax formula.
\end{Proposition}
In the reference, the convergence is stated to be point-wise which is not the case even for $t=0$. Nonetheless, the proof yields the convergence at continuity points. \\
On the other hand,  Barron and Jensen proved the stability of their solutions in \cite[Corollary 3.9]{BarronJensen}. The proof is provided for a bounded Hamiltonian but can easily be extended to our framework. Thus, we have the following corollary:
\begin{corollary}
    \label{coroll: HJ exist unicite lax}
    Given $u_ 0 \in\mathcal{M}_{+}(\R) $ with finite mass, the unique Barron--Jensen solution $U$ is the function given by the Hopf--Lax formula (\ref{eq: hopf lax formula cvx}) where  
$$
U( 0,x)=\itg_ {(- \infty, x)} u_0(dy)\quad \text { on } \mathbb{R}.
$$
\end{corollary}

Since in our case $\Phi^\prime$ is non-negative, the characteristics always point towards the right direction. The uniqueness result can be generalized by specifying the influence domain.
\begin{Lemma}
    \label{prop: domaine d'influence}
    Assume that $U_0$, $V_0$ are l.s.c bounded functions and that there exists $C> 0$ such that: 
    $$
    U_0(x) = V_0(x), \quad \forall x \in (-\infty, C]
    $$
    then, $U(t, x) =V(t,x)$  for all  $ x < C +t \Phi^\prime(0)  $.
\end{Lemma}
\begin{proof}
    For any $t>0, x \in \R$, we observe that if $y$ is such that $\frac{x-y}{t} \geq \Phi^\prime(0)$ then $\Psi\left(\frac{x-y}{t}\right)=0$. So, the infimum can be taken on the smaller set $y \in (-\infty, x- \Phi^\prime(0)t]$. This implies the lemma.
\end{proof}
\begin{remark}
    We emphasize that all these results can be adapted in the case of a strictly concave flux by replacing l.s.c with u.s.c and maximum by minimum. It is sufficient to consider $V(t,x)= U_0(\infty) - U(t,-x).$ The Hopf-Lax formula becomes:
    \begin{equation}
        \label{eq: hopf lax formula concave}
        U(t,x)=\sup _{y \in \mathbb{R}}\left(\itg_{(-\infty,y]} u_{0}(dz) - t \Psi\left(\frac{x-y}{t}\right)\right),
    \end{equation}
where $\Psi$ is the convex conjugate of $-\Phi$.
\end{remark}
\subsection{Hydrodynamic limit on $\Z$}
\label{subsec: hydro limit on Z}
In this subsection, we consider the \ZRP on $\Z$. Let $\mathcal{M}_{+}(\R)$ be the space of non--negative measures on $\R$ endowed with the weak topology.  \\
We consider $(\eta_{N,0})_N \in \N^\Z$ a sequence of configurations  and we define the empirical measure associated with the process $\eta_{N,t}$:
$$ \rho_{N,t} (dy) = \frac{1}{N} \sum_{i \in \Z}  \eta_{N,t}(i) \delta_{i/N}(dy) .$$
Note that we extended the empirical measure to the hole line $\Z$.


\begin{theorem}[Hydrodynamic limit on $\Z$]
    \label{thm: hydrodynamic limit on Z}
    Given $u_0 \in \mathcal{M}_{+}(\R)$ with a finite mass, we consider the Barron--Jensen viscosity solution to  
    \begin{equation}
    \left\{\begin{array}{l}
    \partt U+ \Phi(\partx U)=0, \quad t,x \in \Omega \\
    U(0,\cdot)=U_0,
    \end{array}\right.
\end{equation}
where $U_0$ is the cumulative distribution function of $u_0$.
    Let $ (\rho_N(0, dy))_N $ be a sequence weakly convergent in law to $u_0(dx)$.\\
    We assume that there exists $A>0$ such that:
    \begin{equation}
        \label{eq: assumption support hydro limit Z}
        \lim _{N \rightarrow \infty} \P{\rho_N(0, \R \setminus [-A,A])=0}=1.
    \end{equation}
    Then, the sequence $(t\rightarrow \rho_{N,\frac{Nt}{p-q}})_N$ converge in law in $\mathbb{D}(\R_+, \mathcal{M}_+(\R))$ to a measure-valued function $t \rightarrow u(t,dx)$.
    More over, for each $t \geq 0 $, the c.d.f of $(u_t(dx))$ is $(U_t)^ \star$, the u.s.c envelope of $U_t$.
\end{theorem}

Rezakhanlou \cite{Rezakhanlou1991HydrodynamicLF} has proved the hydrodynamic limit  when the initial configuration is sampled according to a local equilibrium 
 $\otimes_{i \in\Z} \nu\left(u_{i, N}\right)$
 where 
$$\lim _{N \rightarrow \infty} \int_{|x| \leq \ell}\left|u_{[x N], N}-u_0(x)\right| d x=0$$
and $u_0$ has a bounded density with respect to the Lebesgue measure. Actually, in \cite{Rezakhanlou1991HydrodynamicLF}, $g$ is assumed to be bounded. We underline that one can adapt the proof  \cite[Theorem 0.3, Chapter 8]{kipnis1999sli} to generalise the result to non-bounded $g$ when $(\rho_{N,0})$ are compactly supported. Our theorem generalizes the hydrodynamic limit to irregular initial data, as we want to consider the macroscopic behaviour for $u_0=\delta_0$. \\

Let us denote  $\mathcal{Q}_N$ the law of $t \mapsto \rho_{N,\frac{Nt}{p-q}}$. The proof is based on a tightness and uniqueness argument. 
\begin{Lemma}[Tightness]
     The sequence $(\mathcal{Q}_N)_N$ is tight. Moreover, if $\mathcal{Q}$ is any limit point of $(\mathcal{Q}_N)_N$, then for $\mathcal{Q}$ almost every $u$, $u$ is weakly continuous in $t$.
\end{Lemma}
The proof of this lemma is a standard argument similar to the tightness proof on the segment (Lemma \ref{Lemma: tension mesure emirique segment}) and we omit it.
\begin{proof}[Proof of Theorem \ref{thm: hydrodynamic limit on Z}]
     Let $\mathcal{Q}^\star$ be a limit point of $\mathcal{Q}_N$. The idea is to trap the process $\eta$ between two processes with more regular initial data. For any $\varepsilon>0$, we introduce  the following functions
\begin{align*}
    \underline{U}_ 0^\varepsilon (x) &= \inf_y\left \{U_0(y)+\frac{|x-\varepsilon-y|^2}{\varepsilon^2} - \varepsilon \right\}_ + ,\\ 
     \bar{U}_ 0^\varepsilon (x) &= \iota_ \varepsilon(x) + \sup_y
     \left \{U_0(y)-\frac{|x+\varepsilon-y|^2}{\varepsilon^2}\right\}.
\end{align*}
where $\iota_\varepsilon$ is a smooth function equal to $0$ on $(-\infty,-A-2]$ and $\varepsilon$ on $[-A-1,\infty) $. We denote their associated unique viscosity solution by $\underline{U}^\varepsilon$ and $\bar{U}^\varepsilon$.   
We also denote by 
$$
 \inf_\varepsilon \bar{U}^\varepsilon = \Bar{U} \qquad \text{et} \qquad \sup_\varepsilon \underline{U}^\varepsilon = \underline{U}
$$

For any $\varepsilon>0$, $\underline{U}_ 0^\varepsilon$ and $\bar{U}_ 0^\varepsilon $ are Lipschitz, non-decreasing, bounded, and converge to zero when $x$ goes to $-\infty$. We can thus consider two \ZRP with the following initial configurations:

$$\bar \xi_0^\varepsilon \sim \otimes_{i \in\Z} \nu\left(N \left[ \bar U_{0,\frac{i}{N}+\frac{1}{2N}}^\varepsilon - \bar U_{0,{\frac{i}{N}-\frac{1}{2N}}}^\varepsilon \right]\right) \quad \text{and} \quad \underline \xi_0^\varepsilon \sim \otimes_{i \in\Z} \nu\left(N \left[ \underline U_{0,\frac{i}{N}+\frac{1}{2N}}^\varepsilon - \underline U_{0,\frac{i}{N}-\frac{1}{2N}}^\varepsilon\right]\right).$$ 

By \cite[Theorem 1.3]{Rezakhanlou1991HydrodynamicLF}, the hydrodynamic limit of $\underline \xi^\varepsilon$ (resp. $\bar \xi^\varepsilon$) corresponds to the solution of \ref{eq: HJ on R} with initial data $\underline U_ 0^\varepsilon$ (resp $\bar U_ 0^\varepsilon$). \\

Let us denote by $\mathcal{E}_N$, the event where for any $x\in \R$,
\begin{gather*}
     U_0\left(x-\frac{\varepsilon}{2}\right) -\frac{\varepsilon}{2} \leq h(Nx) \leq U_0\left(x+ \frac{\varepsilon}{2}\right) +\frac{\varepsilon}{2}  ; \\
     \bar h _0^\varepsilon (Nx) \geq \bar{U}_0^\varepsilon\left(x-\frac{\varepsilon}{2}\right) -\frac{\varepsilon}{2} ;\\
     \underline h_0^\varepsilon(Nx) \leq \underline U_ 0\left(x+\frac{\varepsilon}{2}\right) +\frac{\varepsilon}{2} .
\end{gather*}
By using the Levy distance, we obtain that $\mathcal{Q}_N(\mathcal{E}_N)$ converges to $1$.
Thanks to the choice of $\bar U_ 0^\varepsilon, \underline U_ 0^\varepsilon$  and the assumption (\ref{eq: assumption support hydro limit Z}), we obtain on the event $\mathcal{E_N}$:
$$
\underline h_0^\varepsilon \leq h_0 \leq \bar h_0^\varepsilon.
$$
A similar coupling as (\ref{eq: coupling attractivity}) ensures for any $t\geq 0 , \varepsilon>0$, and almost surely on $\mathcal{E}_N$,
$$
 \underline h_t^\varepsilon \leq h_t \leq \bar h_t^\varepsilon.
$$
By using the hydrodynamic limit of $\underline \xi^\varepsilon$ and $\bar \xi^\varepsilon$, we obtain for any $t\geq 0$, 
\begin{equation}
     Q^\star \text{-a.s},\quad \forall x \in \operatorname{Cont}(U_t), \qquad \underline U^\varepsilon _t \leq U_t \leq \bar U ^\varepsilon_t
\end{equation}
where $\operatorname{Cont}(U_t)$ are the continuity point of $U_t$. By making $\varepsilon$ go to $0$, we have for any $t\geq 0$:
\begin{equation}
     Q^\star \text{-a.s},\qquad  \forall x \in \operatorname{Cont}(U_t), \qquad \underline U_ t \leq U_t \leq \bar U_t.
\end{equation}

By Proposition \ref{prop cvg formule de LAX},  $\underline{U}(t,x) = \bar{U}_\star(t,x)$ on their continuity point. Thus, we obtain  $\forall t \geq 0,$
$$ Q^\star \text{-a.s},\qquad     \forall x \in \operatorname{Cont}(U_t), \qquad  U_t = \underline U_t  $$
which uniquely determines $\mathcal{Q}^\star$.
\end{proof}
\section{Hydrodynamic limit on the segment}
\label{sec: hydro limit on the segment}
Let $(\eta_{N,0}) \in \Omega_{N,k}$ such that the associated empirical measure $(\rho_{N,0})$ weakly converge to $u_0 \in \mathcal{M}_+([0,1])$. The proof is also a tightness--uniqueness argument. We denote by $\mathcal{Q}_N$ the law of $t \rightarrow \rho_{N,\frac{N}{p-q}t}$.
\begin{Lemma}[Tightness]
    \label{Lemma: tension mesure emirique segment}
     The sequence $(\mathcal{Q}_N)_N$ is tight. Moreover, if $\mathcal{Q}$ is any limit point of $(\mathcal{Q}_N)_N$, then for $\mathcal{Q}$--almost every $u$, $u$ is weakly continuous in t.
\end{Lemma}
The main argument is the point-wise convergence of the cumulative distribution functions.
\begin{Proposition}[Concentration on BJ solution] 
    \label{prop uniqueness cdf on segment}
    Given $(\eta_{N,0}) \in \Omega_{N,k}$ such that the associated empirical measure $\rho_{N,0}$ weakly converge to $u_0 \in \mathcal{M}_+([0,1])$, we consider $U$ the unique Barron--Jensen viscosity solution associated to (\ref{eq: HJ on R}) where $U_0$ is the c.d.f of $u_0$.\\
    For any $t\geq 0, \delta >0$ and  $ y  \in (0, 1) \cap \operatorname{Cont}(U_t)$, we obtain  
    $$\lim _{\substack{N \rightarrow \infty}} \mathbb{P}\left( \left | \frac{1}{N}{h}_{N,\frac{Nt}{p-q}}(y N) - U(t, y ) \right |  \geq \delta  \right)=1 .$$
\end{Proposition}
The proof of Lemma \ref{Lemma: tension mesure emirique segment} and Proposition \ref{prop uniqueness cdf on segment} are deferred to the next subsections. 
\begin{proof}[Proof of Theorem \ref{thm hydro limit on the segment}]
    Given $\mathcal{Q}^\star$ a limit point of $\mathcal{Q}^N$, the previous lemma implies that for any $t\geq 0$,
        $$
    \mathcal{Q}^\star- a.s,     \qquad y \in (0,1) \cap \operatorname{Cont}(U_t), \, \itg_{- \infty} ^ y \pi_{t}^\star(dy) = U(t,y)
    $$
    The total number of particles is conserved. Thus, for any $t\geq 0$,
        $$
    \mathcal{Q}^\star- a.s, \quad \pi_{t}^\star([0,1])= \alpha
    $$
    which uniquely determine $\mathcal{Q}^\star$.
\end{proof}
\subsection{Proof of Lemma \ref{Lemma: tension mesure emirique segment}:}
   
 By \cite[Proposition 1.7]{kipnis1999sli}, it is sufficient to show that  for any  $J \in \mathcal{C}^\infty([0,1])$, the sequence $(t \rightarrow \langle  \rho_{N,\frac{N}{p-q}t},J \rangle $ is tight. Since the total mass of the empirical measure is bounded, it remains to prove that for any $T>0$:
    \begin{equation}
        \label{eq: control tightness 0}
        \varlimsup_{\delta \downarrow 0} \varlimsup_{N \rightarrow \infty} \mathbb{E}\left[\sup _{s, t \leqslant T,|t-s| \leqslant \delta}\left|\left\langle \rho_{N,\frac{N}{p-q}t}-\rho_{N,\frac{N}{p-q}s}, J \right\rangle\right|\right]= 0
    \end{equation}
Let $f(\eta) = \frac{1}{N} \sum_{x=1}^N J\left(\frac{x}{N}\right) \eta(x)$. We consider the associated martingale $(M_{N,t})$ such that for any $t\geq 0$
 \begin{equation}
     f\left ( \eta_{Nt,\frac{Nt}{p-q}} \right) -f\left ( \eta_{N,\frac{Ns}{p-q}} \right) = \itg_{\frac{Ns}{p-q}}^{\frac{Nt}{p-q}} \mathcal{A}_N f\left ( \eta_{N,r} \right) dr + M_{N,\frac{N}{p-q}t}- M_{N,\frac{N}{p-q}s}.
 \end{equation}
A direct computation shows that 

\begin{align}
     &\left|\itg_{\frac{Ns}{p-q}}^{\frac{Nt}{p-q}} \mathcal{A}_N f\left ( \eta_{N,r} \right) dr \right | \notag\\
     &\qquad = \left|\itg_{\frac{N}{p-q}t}^{\frac{Nt}{p-q}} \frac{1}{N} \sum_{x=1}^{N-1} (p g(\eta_r(x))-q  p g(\eta_r(x+1))) \left( J\left(\frac{x+1}{N}\right) - J\left(\frac{x}{N}\right)\right)  dr \right | \notag \\
     &\qquad \leq \|\nabla J\|_\infty g^\star \delta \frac{k}{N}. \label{eq: control tightness 1}
\end{align}

On the other hand, we obtain by Doob inequality:
\begin{align}
    \E{\sup_{s \leq T} M_{N,\frac{Ns}{p-q}}^2}
    & \leq \E{ M_{N,\frac{NT }{p-q}}^2} \notag
    \\&= \E{\itg_0^{\frac{NT}{p-q}} \mathcal{A}_Nf\left ( \eta_{N,r} \right)^2 -2 f\left ( \eta_{N,r} \right) \mathcal{A}_N f\left ( \eta_{N,r} \right) } \notag
    \\&=\frac{1}{N^2 } \E{\itg_0^{\frac{NT}{p-q}} \sum_{x=1}^{N-1} (pg(\eta_r(x)) +q g(\eta_r(x+1))) \left( J\left(\frac{x+1}{N}\right) - J\left(\frac{x}{N}\right)\right)^2 } \notag
    \\ &\leq \frac{k C_{J,T,g}}{(p-q)N^3}.  \label{eq: control tightness 2}
\end{align}
Clearly, (\ref{eq: control tightness 1}) and (\ref{eq: control tightness 2}) imply (\ref{eq: control tightness 0}). \qed
\subsection{Proof of Proposition \ref{prop uniqueness cdf on segment} (Lower bound)}

The idea is to compare the process on the segment to the process on $\Z$ through two different couplings. \\

Given $t^\star,\delta,\varepsilon > 0$, we fix $R= \frac{6 g^\star t^\star}{p-q}$. We denote by $C$ the constant given by Lemma \ref{lem local control in coupling ray - line 1} for $a=0$ and $b=+\infty$. As the macroscopic comparisons are valid only over a period of $\frac{N C \varepsilon}{ (p-q)}$, we introduce a time subdivision $(t_i)_{i\in \N} = (\frac{ C \varepsilon }{ 2} i)_{i\in \N} $. For convenience, let us denote $(T_i)_{i\in \N} = (\frac{ C N \varepsilon }{ 2(p-q)}  i)_{i\in \N} $. \\

We introduce  a first process $\Tilde{\eta}$ on the subgraph $\N$ started at $\Tilde{\eta}_0= \eta( \cdot- \lfloor 2\varepsilon N\rfloor )$ which we couple with the ZRP process $\eta$ on the segment following (\ref{eq: coupling attractivity}). For any $t \geq 0$, 
\begin{equation}
\label{eq: couplage lower bound prop limit hydro}
\Tilde{h}_{N,t} \leq {h}_{N,t}, \qquad \text{a.s}.
\end{equation}

Then, we define $\hat \eta $ a \ZRP on $\Z$ coupled with $\tilde{\eta}$ according to (\ref{eq: coupling attractivity}) on $[T_i, T_{i+1})$ and updated at each $T_i$ by imposing $\hat{\eta}= \tilde \eta$. We claim that for any $t \leq 2t^\star$, the sequence $(\hat \pi(t,dy))_N$ converge in law for the weak topology to the measure $\partx [U_t^\star(\cdot -2\varepsilon)]$.\\

Thanks to Theorem \ref{thm: hydrodynamic limit on Z}, it is sufficient to prove the claim for each update $t_i \leq 2 t^\star $. We argue by induction and we assume that the claim holds at $t_{i-1}$.\\
Let us denote $\mathcal{E}=\{\forall s\leq t^\star, \, \Tilde{h}(s, RN) = k \}$. Thanks to Lemma \ref{lem:attractivity}, we can couple $\tilde \eta$ with $k$ independent particles jumping only to the right with rate $pg^\star$. Thus, 
$$
 \P{\mathcal{E}}\geq 1- \P{\max_{i=1,\cdots,k} X_i \geq RN},
$$
where $X_i$ are independent Poisson variables with parameters $\frac{2 g^\star N t^\star}{p-q}$.
Lemma \ref{lem concentration poisson} implies that r.h.s term converges to $1$.\\
Given a bounded continuous test function $J$, we obtain on the event $\mathcal{E}$ by the triangular inequality: 
    \begin{align}
\left| \frac{1}{N} \sum_{x\in\Z} J\left (\frac{x}{N}\right ) \hat{\eta}\left( [T_i]^+  ,x\right) - \itg J(x+2\varepsilon) \partx U_{t_i}^\star \left(dx\right)\right | \leq A_ {i,N}+B_ {i,N}+C_ {i,N}+D_ {i,N},
\end{align}
where \begin{align*}
     A_{i,N}=& \frac{ 2 \|J\|_\infty}{N} \sum_{x =\varepsilon N } ^{N R }\left|\tilde \eta \left( T_ i  ,x\right)-\hat{\eta}\left( \left[T_ i \right]^-,x\right)\right|\\
     B_{i,N}=& \frac{k}{N} \sup_{x \leq \varepsilon} \left|J(x) -J(0)\right| \\
     C_{i,N}=& \left| \frac{1}{N} \sum_{x= \varepsilon N}^{R N} J(x/N) \hat{\eta}\left( \left[T_ i \right]^-,x\right) - \itg_{-\varepsilon}^{R-2\varepsilon}J(x+2\varepsilon) \partx U_{t_i}^\star \left(dx\right)\right |\\
     D_{i,N}=& J(0) \left| \frac{1}{N}\sum_{x=\varepsilon N}^{R N}  \hat{\eta}\left( \left[T_ i \right]^-,x\right) - \frac{k}{N}\right | .
\end{align*}
First, $(A_{i,N})$ converges to $0$ in expectation by Lemma \ref{lem local control in coupling ray - line 1}. The second term vanishes by continuity of $J$ in $0$. By the induction assumption, $(C_{i,N})$ and $(D_{i,N})$ converge in probability to zero which concludes the proof of the claim. \\

Consequently, for any $t \leq 2t^\star$ and $x \in \R$ such that $x-2\varepsilon \in \operatorname{Cont}(U_t)$:
\begin{align*}
    &\P{\left| \Tilde{h}_{\frac{Nt}{p-q}}(Nx) -U(t,x-2\varepsilon)\right| \geq \frac{\delta}{2}  }
    \leq \P{\mathcal{E}^c}+ \frac{16}{\varepsilon} \mathbb{E}\left[\frac{1}{N} \sum_{x \in\llbracket \varepsilon N, R N \rrbracket}|\tilde \eta_{\frac{Nt}{p-q}}( x)-\hat{\eta}_{\frac{Nt}{p-q}}( x)|\right]  \\
    & \qquad + \P{\left|\sum_{j= \varepsilon N}^{RN}\hat \eta_{\frac{Nt}{p-q}}(j) - \alpha \right|  \geq  \frac{\delta}{8}}+ \P{ \left |\hat{h}_{\frac{Nt}{p-q}}(Nx) - U(t,x-2\varepsilon) \right |\geq \frac{\delta}{8}  } .
\end{align*}
where $\mathcal{E}^c$ the complementary event.
The third and fourth right-hand side terms converge to zero thanks to the hydrodynamic limit of $\hat \eta$. Lemma \ref{lem local control in coupling ray - line 1} guarantees the convergence of the second term.  
Thus, 
$$\P{\left| \Tilde{h}_{\frac{Nt}{p-q}}(Nx) -U(t,x-2\varepsilon)\right| \geq \frac{\delta}{2}  }
     \xrightarrow[N \rightarrow \infty]{} 0.$$
By (\ref{eq: couplage lower bound prop limit hydro}), we obtain
$$\lim _{\substack{N \rightarrow \infty }} \mathbb{P}\left( \frac{1}{N}{h}_{N,t}(x N) \geq U(t, x-2 \varepsilon )  - \frac{\delta}{2}  \right)=1 .$$
Finally, $\operatorname{Cont}(U_t)$ is dense and for any $x \in \operatorname{Cont}(U_t)$, we can  choose   $\varepsilon >0$ such that $$x- 2 \varepsilon \in \operatorname{Cont}(U_t) \qquad \text{and} \qquad U(t,x-2\varepsilon) \geq U(t,x)-\frac{\delta}{2}.$$ Thus,
$$\lim _{\substack{N \rightarrow \infty }} \mathbb{P}\left( \frac{1}{N}{h}_{N,t}(x N) \geq U(t, x)  -\delta  \right)=1, \qquad \forall \delta >0 .$$
\qed

\subsection{Proof of Proposition \ref{prop uniqueness cdf on segment} (Upper bound) }
The proof follows a similar approach to the lower bound proof, but we need to be cautious at the boundary at $N$.
Thanks to Lemma \ref{lem:attractivity}, it is sufficient to prove the same result on the subgraph $\llbracket - \infty , N \rrbracket$. 
\begin{Lemma}
    \label{lem: upper bound on the ray _ hydro limit}
    Let $\Tilde{\eta}_{N}$ a \ZRP on $\llbracket - \infty , N \rrbracket$ started at $\eta_{N,0}$.
    For any $t, \delta>0$ and  $ y  \in (0, 1) \cap \operatorname{Cont}(U_t)$: 
    $$\lim _{\substack{N \rightarrow \infty \\ k / N \rightarrow \alpha}} \mathbb{P}\left( \frac{1}{N}\Tilde{h}_{N,\frac{Nt}{p-q}}(y N) \leq U(t, y )  +\delta  \right)=1 .$$
\end{Lemma}

\begin{proof}[Proof of Lemma \ref{lem: upper bound on the ray _ hydro limit}]
    Given $\delta, t^\star >0$ and $\varepsilon>0$. We fix $R=\frac{6 g^\star  t^\star}{p-q}$ and $C$ the constant provided by the Lemma \ref{lem local control in coupling ray - line 1} for $a=-\infty $ and $b=1$. As above, we introduce a time subdivision $(t_i)_{i\in \N} = (\frac{ C \varepsilon }{ 2} i)_{i\in \N} $ and $(T_i)_{i\in \N} = (\frac{ C N \varepsilon }{ 2(p-q)}  i)_{i\in \N} $.\\

    Whenever we use Lemma \ref{lem local control in coupling ray - line 1}, we lose information about the density profile between $(1-\varepsilon)N$ and $N$. Unlike the lower bound, this part may contain a non-negligible number of particles. To get around this, we introduce a delayed version at each $t_i$ where we gather the particles at a well-chosen point before $N (1-\varepsilon)$. This is still sufficient as we are looking for an upper bound. More precisely, let us define the function $V^d$ by:
 \begin{itemize}
     \item  For any $t \in [t_0,t_1)$, $V_t^d=U_t$ where $U_t$ is the Barron--Jensen viscosity solution to (\ref{eq: HJ on R}) with $U_0$ as initial condition.
     \item For any $i \in \N$, there exists $y_i \in \operatorname{Cont}{V^d(t_i^-)} \cap \left( 1-2 \varepsilon, 1- \varepsilon\right)$
         \[
\left \{
\begin{array}{lr}
    V^d(x,t_i)=  V^d(x,t_i^-) & \text{if} \quad x \leq y_ i, \\
   V^d(x,t_i)= \alpha & \text{if} \quad x > y_i. \\ 
\end{array}
\right.
\]
For any $t \in (t_i , t_{i+1}]$, $V^d$ is equal to the Barron--Jensen viscosity solution started at $V^d(t_i)$.
 \end{itemize}
As the updates at each $t_i$ affect only the function on a support $[1-2\varepsilon, +\infty)$, Lemma \ref{prop: domaine d'influence} implies that $V^d=U$ for any $t\geq 0$ and $x \leq 1-4 \varepsilon$.\\

We also define the microscopic analog $\tilde{\eta}^d$ to $V^d$ on $\rrbracket - \infty, N \llbracket $

\begin{itemize}
    \item For any $t \in [0,T_1)$, $\tilde \eta^d $ is set to be equal to $\tilde \eta$.
    \item For any $i \in \N$, conditionally on $\sigma\left( \tilde \eta(s) , s \leq T_i\right)$
    \[
\left \{
\begin{array}{lr}
   \tilde \eta^d\left(x,T_i\right)= \tilde \eta^d(x,T_i^-) & \text{if} \quad x < \lfloor N y_i \rfloor,\\
   \tilde \eta^d(x, T_i) = \sum_{y= \lfloor N y_i \rfloor}^N \tilde \eta^d(y ,T_i ^-) &  \text{if} \quad x = \lfloor N y_i \rfloor,\\
   \tilde \eta^d(x,T_i)=  0  &  \text{if} \quad x > \lfloor N y_i \rfloor .
\end{array}
\right.
\]
This modification delays the evolution. More precisely, $\tilde h^d(T_i^+,\cdot) \geq \tilde h(T_i^+,\cdot) $. For any $t \in (T_i, T_{i+1}]$, both processes evolve according to (\ref{eq: coupling attractivity}) to respect the order of height functions.
\end{itemize}
We finally introduce a \ZRP $\hat{\eta}^d$ evolving on $\Z$ coupled with $\Tilde{\eta}^d$ according to \ref{eq: coupling attractivity} and updated on each $T_i$ as for $\hat{\eta}^d$. We claim that 
\begin{Lemma}
    \label{lem: convergence de la delayed on Z 1}
    For any $t \leq 2 t^\star$, the sequence $\left(\hat \pi^d\left(\frac{Nt}{p-q},dy\right)\right)_N$ converges in law (in the weak topology) to the measure $\partx V_t^\star$.
\end{Lemma}
We defer the proof of this lemma to the end of this sub-section.
 For any $t\leq 2t^\star$, $x \in \Cont{U_t} \cap (-\infty, 1-4\varepsilon)$:
 \begin{align*}
    &\P{\left| \Tilde{h}^d_{\frac{Nt}{p-q}}\left(Nx\right) -U(t,x)\right| \geq \delta  }
    \leq \P{\Tilde{h}_{\frac{Nt}{p-q}}^d(-R N)  \geq \frac{\delta N}{3}}  \\
    & \quad + \frac{3}{\delta} \mathbb{E}\left[\frac{1}{N} \sum_{j \in\llbracket -R N , N (1-\varepsilon) \llbracket}|\tilde \eta^d_{\frac{Nt}{p-q}}( j)-\hat{\eta}^d_{\frac{Nt}{p-q}}( j)|\right] + \P{ \left |\hat{\pi}^d_{\frac{Nt}{p-q}}(\llbracket-NR,Nx\rrbracket) -U(t,x) \right |\geq \frac{\delta}{3}  }.
\end{align*}
The second right-hand side term vanishes as $N$ goes to infinity thanks to Lemma \ref{lem local control in coupling ray - line 1} and Lemma \ref{lem: convergence de la delayed on Z 1} implies the convergence to zero of the third term. \\
By Lemma \ref{lem:attractivity}, we can couple the process $\tilde \eta^d$ with a process of $k$ independent particles that jumps only to the left with a rate $q g^\star $. Thanks to this coupling, we obtain for any $s\leq \frac{2  t^\star}{p-q}$,
\begin{equation}
    \label{eq concentration_1}
    \P{ \tilde h^d_ {sN} (-RN) \geq 1} \leq \P{ \max_{i=1,\cdot,k} X_i \geq RN },
\end{equation}
where $(X_i)_{=1}^k$ are independent Poisson variables with parameter $\frac{2 q g^\star N t^\star}{p-q}$.
Lemma \ref{lem concentration poisson} implies for any $t\geq 0$ and  $x \in \Cont{U_t} \cap (-\infty, 1-4\varepsilon)$:
 \begin{align*}
    \lim _{\substack{N \rightarrow \infty }}  \P{\left| \Tilde{h}^d_{\frac{Nt}{p-q}}\left(Nx\right) -U(t,x)\right| \geq \delta  } =0.
    \end{align*}
Since $\tilde h^d \geq \tilde h$, we
conclude the proof, up to Lemma  \ref{lem: convergence de la delayed on Z 1}.
\end{proof}

We turn to prove the last lemma. 
\begin{proof}[Proof of Lemma \ref{lem: convergence de la delayed on Z 1}]
    We first recall that for any $i \leq i ^\star$, $\{ \tilde h^d_ {T_i^+} (-RN) \geq 1\}$ occurs with probability converging to one. Thanks to Theorem \ref{thm: hydrodynamic limit on Z}, it remains to prove the convergence for each update $t_i \leq 2 t^\star $.\\
    
    As above, we argue by induction and we assume that the lemma holds at $t_{i-1}$.  Given a bounded continuous test function $J$, we obtain almost surely on  $\{\tilde h^d_ {T_i^+} (-RN) \geq -RN\}$:
    \begin{align}
\left| \frac{1}{N} \sum_{x\in\Z} J\left (\frac{x}{N}\right ) \hat{\eta}^d\left( T_i  ,x\right) - \itg J(x) \partx V_{t_i}^\star \left(dx\right)\right | \leq A_ {i,N}+B_ {i,N}+C_ {i,N}+D_ {i,N},
\end{align}
where \begin{align*}
     A_{i,N}=& \frac{ 2 \|J\|_\infty}{N} \sum_{x \in\llbracket -R N , N (1-\varepsilon) \rrbracket}\left|\tilde \eta^d\left( \left[T_ i \right]^-,x\right)-\hat{\eta}^d\left( \left[T_ i \right]^-,x\right)\right|\\
     B_{i,N}=& \left| J\left(\frac{\lfloor N y_i \rfloor}{N} \right) \frac{k}{N}- J(y_i) \alpha \right  |\\
     C_{i,N}=& \left| \frac{1}{N} \sum_{x=-N R}^{y_i N} J(x/N) \hat{\eta}^d\left( \left[T_ i \right]^-,x\right) - \itg_{-R}^{y_i}J(x) \partx V_{t_i^-}^\star\left(dx\right)\right |\\
     D_{i,N}=& \left| \frac{J\left(\frac{\lfloor Ny_i \rfloor}{N} \right)}{N} \sum_{x=-N R}^{y_i N}  \hat{\eta}^d\left( \left[T_ i \right]^-,x\right) - J(y_ i) \itg_{-R}^{y_i}  \partx V_{t_i^-}^\star\left(dx\right)\right | .
\end{align*}
First, $(A_{i, N})$ converges to $0$ in expectation by Lemma \ref{lem local control in coupling ray - line 1}. The second term converges to $0$ by continuity of $J$. By the induction assumption and Theorem \ref{thm: hydrodynamic limit on Z}, we obtain the convergence in law of $\left(\hat \pi^d\left(T_i^-,dy\right)\right)_N$ to $\partx V_{t_i^-}$. Thus, $(C_{i,N})$ and $(D_{i,N})$ converge to $0$ in probability as $V_{t_i^-}$ is continuous in $y_i$ and $-R$.
\end{proof}

\section{Proof of Lower and Upper bounds on the mixing time } 
\label{sec: mixing time proof}
In this section, we prove the main theorems on the mixing time.  \\
\indent Let us recall that $\wedge_{N,k}= k \delta_1$ and $\vee_{N,k}=k \delta_N$. For the sake of clarity, we will omit the indices when it is clear from the context. As mentioned, $\wedge$ corresponds to the configuration farthest from the equilibrium. By Theorem \ref{thm hydro limit on the segment}, its macroscopic density evolution is described by $U^\alpha$ the unique Barron--Jensen solution of (\ref{eq: HJ on R})  and its macroscopic equilibrium time  $T_{eq,p,\alpha}$ corresponds to the mixing time. 
\subsection{Lower bound of the mixing time}
The proof of the lower bound relies on the hydrodynamic limit of the \ZRP on the segment and uses the position of the left-most particle as a distinctive statistic.     \\
Given a configuration $\eta \in \Omega_{N,k}^0$, we denote by $\ell(\eta)$ the position of the left-most particle. The following lemma controls $\ell(\eta)$ when $\eta$ is sampled according to the equilibrium measure $\pi_{N,k}$. The proof is the analogue of \cite[Lemma 11]{Labb__2019} and  \cite[Proposition 11]{Levin_2016}.
\begin{Lemma}[Concentration of $\ell(\eta)$]
    \label{Lemma equilibriem concentration }
    There exists $C_g$ depending on $g,p$ such that for any $N,k, \delta \in \N$:
    \begin{equation}
        \pi_{N,k}\left( \ell(\eta) \leq N - \delta \right) \leq C_g g(k) \left(\frac{p}{q} \right )^{-\delta}.
    \end{equation}
\end{Lemma}
\begin{proof}
    We observe that 
    $$
\frac{\pi_{N,k}(\eta^{i,N})}{\pi_{N,k}(\eta)}= \left(\frac{p}{q}\right)^{N-i} \frac{g(\eta(i))}{g(\eta(N)+1)} \geq \frac{g(1)}{g(k)}  \left(\frac{p}{q}\right)^{N-i},
$$
where $\eta^{i,N}$ is the configuration obtained by moving a particle from a position $i$ is moved to the position $N$.
Thus, 
\begin{align*}
    \pi_{N,k}(\ell(\eta) =i) =  \sum_{\eta, \ell(\eta)=i}\pi_{N,k}(\eta) \leq \frac{g(k)}{g(1)}  \left(\frac{p}{q}\right)^{-(N-i)} \sum_{\eta, \ell(\eta)=i}\pi_{N,k}(\eta^{i,N}) \leq \frac{g(k)}{g(1)}   \left(\frac{p}{q}\right)^{-(N-i)}.
\end{align*}
We obtain the lemma by summing over $i \leq N-\delta$. 
\end{proof}
We can now prove Theorem \ref{thm lower bound mixing time }.
\begin{proof}[Proof of Theorem \ref{thm lower bound mixing time }]
    As mentioned above, it is natural to consider as initial condition $\wedge$. We obtain for any $\delta>0$,
    \begin{align*}
d_{N, k}\left((T_{eq,p,\alpha}-\delta)N\right) &\geq \VT{ P^\wedge_{(T_{eq,p,\alpha}-\delta)N} - \pi_{N,k}} \\
&\geq \mathbb{P}^{\wedge}_{(T_{eq,p,\alpha}-\delta)N}\left( \ell\left({\eta}\right)\leq N-\varepsilon N\right) - \pi_{N, k}\left(\ell (\eta)\geq N-\varepsilon N\right).
    \end{align*}
    where $P_t^\wedge$ is the law of the process started at $\wedge$.
     Since \begin{equation}
    T_{eq,p,\alpha}= \frac{1}{p-q} \inf\{t >0, U^\alpha(t,dx) \leq \alpha \ind_{x>1}\},
\end{equation} the hydrodynamic limit on the segment (Theorem \ref{thm hydro limit on the segment}) implies that for any $\delta>0$, there exists $\varepsilon \in(0,1)$,
     $$\lim _{\substack{N \rightarrow \infty \\ k / N \rightarrow \alpha}} \mathbb{P}^{\wedge}_{(T_{eq,p,\alpha}-\delta)N}\left( \ell\left({\eta}\right)\leq N-\varepsilon N\right)=1 .$$
 Lemma \ref{Lemma equilibriem concentration } concludes the proof.
\end{proof}
\begin{remark}
    Let us emphasise that the same method would apply for $g$ of polynomial growth if the hydrodynamic limit holds.
\end{remark}

\subsection{Upper bound: totally asymmetric jumps and convex flux}
Let us first observe that as $p=1$, the system is mixed when it reaches the configuration $\vee$. Thanks to the attractivity, we have in particular
\begin{align*}
    \forall t \geq 0, \; d_{N,k}(t)= P^\wedge_t\left(h(N-1)>0\right).
\end{align*}
On the other hand, the Hopf--Lax formula written (\ref{eq: hopf lax formula cvx}) as:
\begin{align*}
    U^\alpha(x,t)= \inf\{\inf_{z \leq \frac{x}{t}}[\alpha + t \Psi(z)] ; \inf_{z \geq \frac{x}{t}}[ t \Psi(z)]\}
\end{align*}
implies that the macroscopic equilibrium is reached at 
$$T_{eq,1,\alpha} =\frac{1}{\Phi^\prime(0)}= \frac{1}{g(1)}.$$
As $T_{eq,1,\alpha} $ depends only on the value of $g$ in $1$, the idea is to couple through (\ref{eq: coupling attractivity}) with the linear rate function, the simplest and slowest convex rate function. 

We consider a \ZRP $(\eta^{\operatorname{ind}})$ started at $\wedge$ and a linear rate function $g_{\operatorname{ind}}:x \rightarrow g(1)x$. Since $g\geq g_{\operatorname{ind}}$, we couple both processes  according to (\ref{eq: coupling attractivity}) to obtain almost surely for any $t>0$:
$$ h^\mathrm{ind}_t \geq h_t.$$
We recall that \ZRP with a linear rate corresponds to independent random walks. Consequently, we obtain for any $t\geq 0$,
\begin{align*}
    d_{N,k}(t)= \P{h_t(N-1)>0} \leq  \P{h^\mathrm{ind}_t(N-1)>0} = 1 - \P{\mathcal{P}(tg(1))\geq N}^k .
\end{align*}
If $t$ is fixed to $N \left(\frac{1}{g(1)} +\varepsilon\right)$, the limit is zero, which yields the first assertion of Theorem \ref{thm upper bounds mixing time }. 
\subsection{Upper-bound: Concave flux}
    In this subsection, we will work within the framework of  Theorem \ref{thm upper bdd concave}. The main idea is to show that the leftmost particle and the number of particles at position $N$ approach their equilibrium state around $N T_{eq,p,\alpha}$. Once there, the system will mix rapidly. 
 Let us denote by:
    \begin{equation}
        L_{N,k}=\ell(\wedge_{t}), \qquad \text{and} \qquad S_{N,k}= \wedge_{t}(N),
    \end{equation}
    where $(\wedge_{t})_t$ is a ZRP$(g,p)$ started from the configuration $\wedge$. We also introduce their macroscopic analogue:
    \begin{equation}
        \forall t\geq 0, \qquad 
    \begin{array}{rl}
     \macroleft(t) =  &  \inf \{x \in  \R, U^\alpha((2p-1)t,x)>0\},\\
      &\\
     \macrostack(t)= &\alpha-U^\alpha((2p-1)t,1^-).
    \end{array}
    \end{equation}
    The macroscopic equilibrium time could be written as:
    $$T_{eq,p,\alpha}=\left(\macroleft\right)^{-1}(1)=\left(\macrostack\right)^{-1}(\alpha).$$
    
    The following lemma ensures that the microscopic behaviour corresponds to the expected macroscopic behaviour.
    \begin{Proposition}[Macroscopic behaviour of $L_{N,k}$ and $S_{N,k}$]
        \label{prop control of L and R}
       For any  $t\leq T_{eq,p,\alpha}$, we have the following convergence in probability:
          \[
          \frac{1}{N}S_ {N,k}(Nt)  \xrightarrow[\substack{N \rightarrow \infty \\ k / N \rightarrow \alpha}]{}  \macrostack (t) \qquad \text{ and } \qquad \frac{1}{N}L_ {N,k}(Nt)   \xrightarrow[\substack{N \rightarrow \infty \\ k / N \rightarrow \alpha}]{}  \mathscr{L}^\alpha_ p (t).
          \] 
    \end{Proposition}
     The proof of this lemma is deferred to Section \ref{sec: control de  L et de S}. As in the case of the exclusion process \cite[Proposition 9]{Labb__2019}, we emphasise that the hydrodynamic limit of $(\wedge_t)$ provides only the upper--bound on $L_{N,k}$ and $S_{N,k}$: For any $\delta>0$
     \[
\forall t \leq T_{eq,p,\alpha},\; \lim_{N\rightarrow \infty} \P{\frac{S_ {N,k}(Nt)}{N} \leq \macrostack(t)+\delta \text{ and } \lim_{N\rightarrow \infty} \frac{L_ {N,k}(Nt)}{N} \leq \macroleft(t)+\delta}=1.
     \]
     
     Let us focus on the proof of the second statement of \ref{thm upper bounds mixing time }.
     \begin{proof}[Proof of Theorem \ref{thm upper bdd concave}]
        By triangular inequality, we have for any $t \geq 0$:
        $$ d_{N,k}(t) \leq 2 \max_{h, h^\prime}\left\|P_{t}^h-P_{t}^{h ^\prime}  \right\|_{\operatorname{TV}}.$$
        The attractivity of the \ZRP permits a first simplification. Indeed, starting for any $h,h^\prime$, we can trap their evolution between $(\wedge_t)$ and $(\vee_t)$ thanks to a similar coupling to (\ref{eq: coupling attractivity}):
        \begin{align*}
            d_{N,k}(t) \leq 2 \max_{h, h^\prime}\left\|P_{t}^h-P_{t}^{h ^\prime}  \right\|_{\operatorname{TV}} &\leq 2 \P{ \wedge_{t} \neq \vee_t}.
        \end{align*}
        Then, once the process $(\wedge_t)$ reaches the minimal state $\vee$, the processes $(\wedge_s)_s$ and $(\vee_s)_s$ become identical:
        $$
        d_{N,k}(t) \leq 2 \P{\forall s \leq t, \, \wedge_{s} \neq \vee}.$$
         For any $\varepsilon>0$, we define
         \[
         \mathcal{E}_{N,k}^\varepsilon= \left\{  \eta \in \Omega_{N,k}^0, \; \ell(\eta) \geq N (1-\varepsilon) \text{ and } \eta(N) \geq (1-\varepsilon ) k \right\}.
         \]
        We also denote by $\wedge^\varepsilon$ the maximal element of the set $ \mathcal{E}_{N,k}^\varepsilon$. For any $\delta>0$, we obtain
        \begin{align*}
            & d_{N,k}(N\left[T_{eq,p,\alpha}+\delta\right])\\ & \quad \leq 2 \P{\forall t \leq N\left[T_{eq,p,\alpha} +\delta\right], \, \wedge_{t} \neq \vee} \\
            & \quad \leq 2\P{\forall t \leq N\left[T_{eq,p,\alpha} +\delta\right], \, \wedge_{t} \neq \vee \, \&\,  \wedge_{N T_{eq,p,\alpha} } \in \mathcal{E}_{N,k}^\varepsilon}  + 2\P{\wedge_{NT_{eq,p,\alpha}} \notin \mathcal{E}_{N,k}^\varepsilon} \\
            &\quad \leq 2\P{\forall t \leq \delta N, \,  \wedge^\varepsilon_{t} \neq \vee} +2\P{\wedge_{NT_{eq,p,\alpha}} \notin \mathcal{E}_{N,k}^\varepsilon},
        \end{align*}
        where $(\wedge_t^\varepsilon)$ is a \ZRP started at $\wedge^\varepsilon$.\\
        Since $T_{eq,p,\alpha}= (\macroleft)^{-1}(1)= (\macrostack)^{-1}(\alpha)$, Proposition \ref{prop control of L and R} implies that 
        $$
        \P{\wedge_{NT_{eq,p,\alpha}} \notin \mathcal{E}_{N,k}^\varepsilon} \xrightarrow[\substack{N \rightarrow \infty \\ k / N \rightarrow \alpha}]{} 0.
        $$
        We compare the system to an exclusion process to control the last term. We denote by $\bar g= \max_k g(k)$.
        Let us introduce two \ZRP $(\tilde \wedge_t^\varepsilon)$ (resp. $(\tilde \vee_t)$) with a modified rate function $\tilde{g}=(pg(1)+q \Bar{g}) \ind_{\cdot\geq 1}$ and a bias $\tilde{p}= \frac{pg(1)}{pg(1)+q \bar g}$ started at $\wedge^\varepsilon$ (resp. $\vee$). We couple these three processes with (\ref{eq: coupling attractivity}):
        $$
         \forall t\geq 0,\qquad h_{\tilde \wedge_t} \geq \max(h_{\wedge_t}, h_{\tilde \vee_t}) \qquad \text{  a.s.}
        $$
        We obtain that:
        \begin{align*}
            \P{\forall t \leq \delta N, \,  \wedge^\varepsilon_{t} \neq \vee} & \leq  \P{\forall t \leq \delta N, \,  \Tilde {\wedge}^\varepsilon_{t} \neq \vee}\\ &\leq \P{\tilde \wedge^\varepsilon_{\frac{\delta}{2} N} \neq \tilde\vee_{\frac{\delta}{2} N }} + \P{\forall t \in \left[ \frac{\delta}{2}N, \delta N\right], \, \tilde \vee_t \neq \vee}.
        \end{align*}
        We recall that Lemma \ref{Lemma: link ZRP-exclusion} implies that the  processes $(E^{-1}(\tilde \vee_t))$ and $(E^{-1}(\tilde \wedge_t^\varepsilon))$ evolve according to an accelerated exclusion process.
        So, we can use the same argument as in \cite[Theorem 2]{Labb__2019} and \cite[Lemma 12]{labbé2018mixingtimecutoffweakly} since we are reduced to an exclusion process. \\
        
        In particular, there exists an increasing function $f_{N,k}$ on $\Omega_{N,k}$ and such that:
         $$\mathbb{E}\left[f_{N, k}\left({h_t^{\prime}}\right)-f_{N, k}\left(h_t\right)\right] \leq e^{-\varrho t}\left(f_{N, k}\left(h^{\prime}\right)-f_{N, k}(h)\right),$$
         where $h,h^\prime$ are \ZRP($\tilde g, \tilde p$) coupled according to (\ref{eq: coupling attractivity}) and $\varrho$ is a constant (see \cite[Equation 45]{labbé2018mixingtimecutoffweakly}). We obtain by Markov's inequality and \cite[Lemma 12]{labbé2018mixingtimecutoffweakly}:
         \begin{align*}
             \P{\tilde \wedge^\varepsilon_{\frac{\delta}{2} N} \neq \tilde \vee_{\frac{\delta}{2} N }} \leq \frac{\E{ f_{N,k}\left(h_{\tilde \wedge^\varepsilon_{\frac{\delta}{2} N}}\right) - f_{N,k}\left(h_{\tilde \vee_{\frac{\delta}{2} N }}\right)} }{\min_{h \leq h^\prime} \left[\left(f_{N, k}\left(h^{\prime}\right)-f_{N, k}(h)\right)\right ]} \leq C N \left[\frac{pg(1)}{q \bar g}\right]^{\varepsilon N} e^{-\frac{\delta}{2}\rho N}.
         \end{align*}
         
         Thus, for a fixed $\delta>0$, we can choose $\varepsilon$ small enough to obtain:
        $$
         \P{\tilde \wedge^\varepsilon_{\frac{\delta}{2} N} \neq \tilde \vee_{\frac{\delta}{2} N }} \xrightarrow[\substack{N \rightarrow \infty \\ k / N \rightarrow \alpha}]{} 0.
         $$
        Finally, the Markov inequality implies:
    $$\P{\forall t \in \left[ \frac{\delta}{2}N, \delta N\right], \, \tilde \vee_t \neq \vee} \leq \frac{2 \mathbb{E}_\vee({\tilde \tau_ \vee})}{\delta N} = \frac{2}{\delta N \tilde \pi_{N,k}(\vee)},$$
        where $\tilde \tau_\vee $ an excursion's duration from the configuration $\vee$ and $\tilde \pi_{N,k}$ the equilibrium measure of \ZRP $(\tilde g,\tilde p)$. It is then sufficient to lower--bound  $\tilde \pi_{N,k}(\vee)$ by a constant.
        For any $\eta \in\Omega_{N,k}^0$, we have
        $$\tilde \pi_{N,k} (\vee) \geq \tilde \pi_{N,k}(\eta)$$ 
        Thus, for any $\Delta>0$,
        \begin{align*}
           \tilde \pi_{N,k}(\vee) \, \#\left\{\eta \in \Omega_{N,k}^0, \, \ell(\eta) + \eta(N) \geq N + k- \Delta \right\} \geq  \tilde \pi_{N,k}\left(\ell(\eta) + \eta(N) \geq N + k- \Delta  \right)
        \end{align*}
        Due to \cite[Lemma 11]{Labb__2019}, we can choose $\Delta$ independent from $N,k$ such that $$
        \tilde \pi_{N,k}\left(\ell(\eta) + \eta(N) \geq N + k- \Delta  \right) \geq \frac{1}{2}.$$
        We conclude by observing that
        $$
        \#\left\{\eta \in \Omega_{N,k}^0, \, \ell(\eta) + \eta(N) \geq N + k- \Delta \right\} \leq \sum_{i=1}^\Delta \binom{\Delta}{i}.
        $$
     \end{proof}

     \begin{remark}[Pre--Cutoff]
         We emphasise that by using the same approach, we obtain the pre--cutoff for asymmetric \ZRP as long as $pg(1)>q\bar g$ and $\Phi$ strictly concave (see \cite[Chapter 18]{LevinPeresWilmer2006} for more background on the Pre--Cutoff). Indeed, we can couple the \ZRP to \ZRP$(\tilde g, \tilde p)$ where $\tilde{g}=(pg(1)+q \Bar{g}) \ind_{\cdot\geq 1}$ and  $\tilde{p}= \frac{pg(1)}{pg(1)+q \bar g}.$ By using Lemma \ref{Lemma: link ZRP-exclusion} and \cite[Proposition 9]{Labb__2019}, we have for any $\varepsilon>0$:
          $$
        \P{\wedge_{N\tilde T_{eq,p,\alpha}} \notin \mathcal{E}_{N,k}^\varepsilon} \xrightarrow[\substack{N \rightarrow \infty \\ k / N \rightarrow \alpha}]{} 0,
        $$
        where $\tilde T_{eq,p,\alpha}$ is associated to \ZRP$(\tilde g, \tilde p)$. Then, the same method yields for any $\theta>0$,
        \[
         T_{eq,p,\alpha} \leq \liminf _{\substack{N \rightarrow \infty \\ k / N \rightarrow \alpha}} \frac{T_{\text {mix }}^{N, k}(\theta)}{N}  \leq \limsup _{\substack{N \rightarrow \infty \\ k / N \rightarrow \alpha}} \frac{T_{\text {mix }}^{N, k}(\theta)}{N} \leq  \tilde T_{eq,p,\alpha}
        \]
     \end{remark}
\section{Evolution of the left--most particle and the stack at N}
\label{sec: control de  L et de S}
In the sequel, we prove Proposition \ref{prop control of L and R} and assume that the flux $\Phi$ is strictly concave and bounded by $\bar g$. We recall that 
    \begin{equation}
        L_{N,k}=\ell(\wedge_{t}) \qquad \text{and} \qquad S_{N,k}= \wedge_{t}(N)
    \end{equation}
    where $(\wedge_{t})_t$ is a ZRP$(g,p)$ started from the configuration $\wedge$.
    
Let us first describe the behaviour of $\macroleft$ and $\macrostack$. Thanks to the Hopf--Lax formula (\ref{eq: hopf lax formula concave}), we can explicitly compute $U^\alpha $.
\begin{equation}
   U^\alpha (t,x)= \left\{
\begin{array}{ll}
    \max\left(0, \alpha - t \Psi\left(-\frac{x}{t}\right)\right)& \text{if } x \leq g(1)t, \\
    \alpha  & \text{if } x > g(1)t.
\end{array}
\right.
\end{equation}
This implies that  $\macroleft (t) =\alpha \mathscr{L} \left(\frac{(2p-1)t}{\alpha}\right)$ and $\macrostack (t) =\mathscr{S} \left((2p-1)t\right)$ where 
\begin{equation}
      \mathscr L (t)= \left\{
\begin{array}{ll}
    0& \text{if } t \leq \frac{1}{\Psi\left(0\right) } ,\\
    -t \Psi^{-1}\left(\frac{1}{t}\right) & \text{otherwise }
\end{array}
\right.
\qquad , \qquad 
      \mathscr S (t)= \left\{
\begin{array}{ll}
     0 & \text{if } t \leq \frac{1}{g(1) } ,\\
    t \Psi\left(-\frac{1}{t}\right) & \text{otherwise}.
\end{array}
\right.
\end{equation}
These functions are continuous, increasing and strictly convex. Indeed, we can easily compute their derivatives: 
\[
\mathscr L^\prime (t)= \frac{\Phi\left(\Psi^\prime \left(\Psi^{-1}(1/t)\right)\right)}{\Psi^\prime \left(\Psi^{-1}(1/t)\right)}, \; \forall t \geq   \frac{1}{\Psi(0)}\qquad \text{,} \qquad \mathscr S ^\prime(t)=\Phi\left(\Psi^\prime \left(\frac{-1}{t}\right)\right),  \; \forall t \geq   \frac{1}{g(1)}.
\]
The proof is now divided into two subsections: \ref{subsec control de L} for the lower--bound on $(L_{N,k})$ and \ref{subsec control de R} for the lower--bound on $(S_{N,k})$. Unlike for the exclusion process, we do not have a particle-hole symmetry, so we need to prove both estimates.

\subsection{Convergence to $\macroleft$}
\label{subsec control de L}

Let us start by noting that the fastest density is the null density. Thus, if we suppose that a few particles are left behind $\lfloor N \macroleft \rfloor$, they will travel at a higher speed and eventually catch up. The proof will rely on this observation. We first prove a partial version of the Proposition.
\begin{Lemma}
    \label{Lemma convergence of L partial}
   Given $\varepsilon>0$, let $t^\star(\varepsilon)=\left(\macroleft\right)^{-1}\left(6 g^\star q \varepsilon \right) $. For any $\kappa \in \N$ such that {$t^\star(\varepsilon)+ (\kappa -1)\varepsilon \leq T_{eq,p,\alpha}$} and any $\gamma >0$,
        \begin{equation}
            \label{eq induction_1_left_most}
            \tag{$H_{\kappa,1}$}
            \lim_{\substack{N \rightarrow \infty \\ k / N \rightarrow \alpha}} \P{ L_ {N,k}[N t^\star(\varepsilon)+ \kappa \varepsilon N]    \geq N\left[\macroleft\left((\kappa-1)\varepsilon +\frac{\alpha }{\lambda_ 2 (2p-1)}\right)-\gamma\right]} =1
        \end{equation}
        where $\lambda_2= \Psi(0).$
\end{Lemma}
\begin{proof}
    Given $\varepsilon>0$, we introduce the associated time subdivision $T_\kappa =N\left[t^\star(\varepsilon) + \kappa \varepsilon\right]$.
    For convenience, we denote $W_{\kappa,\varepsilon}=\macroleft\left(\kappa\varepsilon +\frac{\alpha }{\lambda_ 2 (p-q)}\right)$ and $R_{\kappa,\varepsilon} = \macroleft(t^\star(\varepsilon)+ \kappa \varepsilon)$. \\
    
    We will prove ($H_{\kappa,1}$) by induction on $\kappa$. 
    Let us sketch the idea of the proof. At  $T_ \kappa$, the hydrodynamic limit implies that there are only $o(N)$ particles between $NW_{\kappa,\varepsilon}$ and $NR_{\kappa,\varepsilon}$. 
    In the totally asymmetric case, we can follow the same approach as in the convex case. Indeed, let us remark that when density converges to zero, the asymptotic speed is $g(1)$. So, we can compare to a \ZRP with a constant jump rate fixed at  $g(1)$ and use the analogous result proven by Labbé \& Lacoin \cite[Proposition 9]{Labb__2019} for \ASEP using the bijection between the two processes. If $p \neq 1$, strict concavity of $\Phi$ allows for an improvement of the result.\\

\begin{figure}[h!]
    \label{fig example tilde eta left most}
    \centering
    \includegraphics{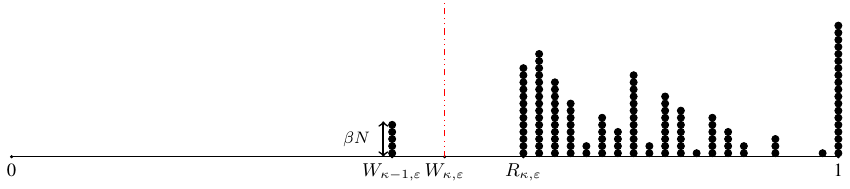}
    \caption{The configuration $\tilde \eta$ at $t=T_\kappa$ zoomed out by a factor $N$. The red dashed line corresponds to the added wall.   }
\end{figure}

    For $\kappa=1$, the result is trivial. Assume ($H_{\kappa,1}$) holds and that {$t^\star(\varepsilon)+ \kappa \varepsilon \leq T_{eq,p,\alpha}$}. Given $\gamma>0$, we introduce three parameters $\gamma_1, \gamma_2, \beta >0$ that we will choose later.\\
    We define a process $\tilde \eta$ (figure \ref{fig example tilde eta left most}) initialized at time $T_ \kappa$ by
        \[\tilde \eta_{T_ \kappa}\left(x\right)=
\left \{
\begin{array}{lr}
    0& \text{if} \quad x \leq  \left  \lceil N R_{\kappa,\varepsilon} - N\gamma_ 2\right \rceil \, \& \, x \neq \left \lceil N W_{\kappa-1,\varepsilon}-N\gamma_1 \right\rceil,\\
   \beta N&  \text{if} \quad x= \left \lceil N W_{\kappa-1,\varepsilon}-N\gamma_1 \right\rceil, \\
   \wedge_ {T_\kappa}(x) &  \text{if} \quad x >\left  \lceil N R_{\kappa,\varepsilon} - N\gamma_ 2\right \rceil.
\end{array}
\right.
\]
Regarding the dynamics, we impose that no particles cross through the site $ \lfloor N W_{\kappa,\varepsilon}\rfloor $. The evolution on the sites after $ \lfloor N W_{\kappa,\varepsilon}\rfloor $ will follow a ZRP$(g,p)$ dynamics. However, the particles on the sites before $ \lfloor N W_{\kappa,\varepsilon}\rfloor $ evolve with a modified rate function $\tilde{g}=(pg(1)+q \Bar{g}) \ind_{\cdot\geq 1}$ and a bias $\tilde{p}= \frac{pg(1)}{pg(1)+q \bar g}$. Let us remark that the generator could be written in the generic form (\ref{eq: generic form generator}). For convenience, $(h_t)$ (resp.$(\tilde h_t)$)  denotes the height function associated to $(\wedge_t)$ (resp. $(\tilde \eta_t)$).
\\
We denote by 
$$\mathcal{E}=\{L_ {N,k}(T_ \kappa)   \geq N\left[W_{\kappa-1,\varepsilon}-\gamma_1\right]\} \cap \{h_{T_ \kappa}\left(\left \lceil N R_{\kappa,\varepsilon} -N \gamma_ 2\right \rceil\right) \leq \beta N\}.$$
Thanks to Theorem \ref{thm hydro limit on the segment} applied to $(\wedge_t)$, we have: $$\lim_{N \rightarrow +\infty}\P{h_{T_ \kappa}\left(\left \lceil N R_{\kappa,\varepsilon} -N \gamma_ 2\right \rceil\right) \leq \beta N}=1. $$
Therefore, the induction assumption $(H_{\kappa,1})$ implies that $\mathcal{E}$ occurs with probability converging to one.\\

On the event $\mathcal{E}$, we observe that the height functions are ordered: 
$$\tilde h_{T_\kappa} \geq h_{T_\kappa}.$$
Since $\tilde{p}\tilde{g} \leq pg $ and $\tilde q  \tilde g= q \bar g \leq qg $,  the processes $(\tilde \eta_ t)$ and $(\wedge_ t)$ coupled through (\ref{eq: coupling attractivity})  verify almost surely on $\mathcal{E}$
    $$
     \forall \tau \geq t \geq T_\kappa, \qquad
    \tilde h_ t \geq  h_t,
    $$
    where $$\tau= \inf\left \{s \geq T_ \kappa, \quad \tilde \eta_s \left (  \left \lceil N W_{\kappa,\varepsilon} \right \rceil +1 \right )\geq 1\right\},$$
 the first time that the $\tilde\eta$-particles  beyond the wall attempt to cross it.\\
 
We obtain the following inequality
\begin{equation}
    \label{eq: control L subdivision evenement}
    \begin{aligned}
    & \P{L_{N,k}\left(T_{\kappa +1}\right)\leq N\left[ W_{\kappa,\varepsilon}-\gamma\right]}\\ & \qquad \leq  \P{ \left\{\ell(\tilde \eta)\left(T_{\kappa +1}\right)\leq N\left[ W_{\kappa,\varepsilon}-\gamma\right]\right\}\cap \mathcal{E}} +\P{{\mathcal{E}^c}} +\P{\tau \leq N\varepsilon}.
\end{aligned}
\end{equation}
Regarding the third term of (\ref{eq: control L subdivision evenement}), we couple through (\ref{eq: coupling attractivity}) to a process where the particles on sites $i>\lceil N R_ {\kappa,\varepsilon} \rceil$ jump to the left independently from each other with rate $qg(1)$. If $\tau \leq  N\varepsilon $, at least one of the particles had jumped more than $\lfloor NR_{\kappa,\varepsilon} - N \gamma_ 2\rfloor- \left \lceil N W_{\kappa,\varepsilon} \right \rceil -1  $. We obtain that
    \begin{align*}
        \P{\tau \leq \varepsilon } \leq \P{\max_{k} X_i \geq \lfloor NR_{\kappa,\varepsilon} - N \gamma_ 2\rfloor- \left \lceil N W_{\kappa,\varepsilon} \right \rceil -1 } ,
    \end{align*}
where $(X_i)_{i=1}^k$ are independent Poisson variable with parameter $qg^\star \varepsilon N$.
As $\macroleft$ is convex,
$$
    R_{\kappa,\varepsilon}- W_{\kappa,\varepsilon} \geq R_{0,\varepsilon}-W_{0,\varepsilon} \geq 6 \varepsilon  g^\star q.
    $$
If we fix $\gamma_2 \leq \varepsilon g^\star q$,  Lemma \ref{lem concentration poisson} implies that $\P{\tau \leq N\varepsilon}$ converge to zero.\\

It remains to control the first term of (\ref{eq: control L subdivision evenement}). The process restricted to  $I_{\left \lceil N W_{\kappa,\varepsilon} \right \rceil}$ behaves like a ZRP$(\tilde{g},\tilde{p})$ on a segment with the following features:

$$ \tilde{N}= \left\lceil N W_{\kappa,\varepsilon} \right \rceil, \qquad \tilde{k} = \beta N,$$

and initialised by $\xi=\beta N \delta_{\left \lceil N W_{\kappa-1,\varepsilon}-N\gamma_1 \right\rceil}$.
Thanks to Lemma \ref{Lemma: link ZRP-exclusion} and  \cite[Proposition $13$]{Labb__2019}, we obtain 
$$
\forall t \geq \frac{(\sqrt{\beta}+\sqrt{W_{\kappa,\varepsilon}-W_{\kappa-1,\varepsilon}+\gamma_1})^2}{pg(1)-q \bar g}, \qquad  \lim_{\substack{N \rightarrow \infty}} \P{ \ell(\tilde \eta)\left(Nt +T_\kappa \right)   \geq N\left[ W_{\kappa,\varepsilon}-\gamma\right]} =1.
$$
To conclude, we should choose $\beta, \gamma_1$ verifying:
\begin{equation}
    \label{eq: critere comparison preuve lefet most}
    \varepsilon \geq \frac{(\sqrt{\beta}+\sqrt{W_{\kappa,\varepsilon}-W_{\kappa-1,\varepsilon}+\gamma_1})^2}{pg(1)-q \bar g}.
\end{equation}
 As $\macroleft$ is convex, we have 
 
 $$ W_{\kappa,\varepsilon} -W_{\kappa-1,\varepsilon} \leq \varepsilon (\macroleft)^\prime (T_{eq,p,\alpha}) =  \varepsilon (p-q)\mathscr{L}^\prime( \mathscr{L}^{-1}(1/\alpha))$$
 and since we assumed $$\frac{pg(1)-q \bar g}{p-q} > \mathscr{L}^\prime( \mathscr{L}^{-1}(1/\alpha)),$$
 we can choose $\beta,\gamma_1$ small enough to verify (\ref{eq: critere comparison preuve lefet most}) which concludes the proof.
\end{proof}
\begin{remark}
    \label{rem: cond  explication 1}
    The condition (\ref{cond sur p et alpha}) is due to the comparison to the exclusion process, which lacks precision. Indeed, we require that the maximal speed attained $(\macroleft)^\prime (T_{eq,p,\alpha})$ is less than the highest speed reachable by \ZRP$(\tilde g,  \tilde p)$. The same method would imply a more general solution if one could control $L_{N,k}$ for small densities. 
\end{remark}
We can now finish the proof of the lower--bound on $L_{N,k}$.
\begin{proof}[Proof of Proposition \ref{prop control of L and R} (Part 1)]
    Given $t \leq T_{eq,p,\alpha}$ and $\varepsilon,\gamma >0$, we fix $$\kappa = \left \lfloor \frac{t- t^\star(\varepsilon)}{\varepsilon} +1\right \rfloor -1 .$$
    As $\macroleft$ is non--decreasing,
    $$
    \macroleft\left((\kappa-1)\varepsilon +\frac{\alpha }{\lambda_ 2 (p-q)}\right) \geq  \macroleft\left (t- t^\star(\varepsilon)-2 \varepsilon +\frac{\alpha }{\lambda_ 2 (p-q)}\right).
    $$
    Since $t^\star(\varepsilon)$ converge to $\frac{\alpha}{\lambda_2(p-q)}$ and $\macroleft$ is continuous, there exist $\varepsilon_ 1(\gamma)$ such that for all $\varepsilon < \varepsilon_1(\gamma)$, $$\macroleft\left((\kappa-1)\varepsilon +\frac{\alpha }{\lambda_ 2 (p-q)}\right) \geq \macroleft(t)- \gamma/3.$$ 
    It follows from the previous lemma that: 
    $$
    \lim_{\substack{N \rightarrow \infty \\ k / N \rightarrow \alpha}} 
    \P{ L_ {N,k}(t^\star(\varepsilon)N+ \kappa \varepsilon N)   
    \geq N\left[\macroleft\left(t\right)-2\gamma/3\right]} =1.$$
    Finally, we observe that
    $$
    t-\varepsilon \leq t^\star(\varepsilon)+ \kappa^\star \varepsilon  \leq t.
    $$
    During $N(t- t^\star(\varepsilon)+ \kappa^\star \varepsilon )$, $L_{N,k}$  cannot move any further than $N\gamma/3$. Indeed, between $t^\star(\varepsilon)N+ \kappa \varepsilon N$ and $t$, we couple $(\wedge_t)$ by (\ref{eq: coupling attractivity}) to the process where particles jump independently to the left with rate $qg(1)$. 
    Thus,
    \begin{align*}
        &\P{ L_ {N,k}(t N) \geq N\left[\macroleft\left(t\right)-\gamma\right]} \\ &\qquad \geq \P{ L_ {N,k}(t^\star(\varepsilon)N+ \kappa \varepsilon N)   \geq N\left[\macroleft\left(t\right)-2\gamma/3\right]}  - \P{\max_{k} X_i \geq \frac{N\gamma}{3}} ,
    \end{align*}
where $(X_i)_{i=1}^k$ are independent Poisson variable with parameter $N( t^\star(\varepsilon)+ \kappa^\star \varepsilon -t)$.
    Lemma \ref{lem concentration poisson} implies the result once $\varepsilon$ is small enough.
\end{proof}

\subsection{Convergence to $ \macrostack$}
\label{subsec control de R}
Let us start with a lemma that guarantees a partial result. 
\begin{Lemma}
    \label{Lemma: convergence of  S _ partial}
     Given $\varepsilon>0$ , let $t^\star(\varepsilon)=\left(\macrostack\right)^{-1}\left(6 \bar g q \varepsilon\right) $. For any $\kappa \in \N$ such that {$t^\star(\varepsilon)+ (\kappa -1)\varepsilon \leq T_{eq,p,\alpha}$} and any $\gamma >0$,
        \begin{equation}
            \label{eq induction_2_left_most}
            \tag{$H_{\kappa,2}$}
            \lim_{\substack{N \rightarrow \infty \\ k / N \rightarrow \alpha}} \P{ S_ {N,k}[N t^\star(\varepsilon)+ \kappa \varepsilon N]    \geq N\left[\macrostack\left((\kappa-1)\varepsilon +\frac{1 }{g(1)(2p-1)}\right)-\gamma\right]} =1
        \end{equation}
\end{Lemma}
Using the same argument as in Lemma \ref{Lemma convergence of L partial} to derive the convergence of $L_{N,k}$, one can show that Lemma \ref{Lemma: convergence of  S _ partial} implies the convergence of $S_{N,k}$ in probability.
\begin{proof}
    Given $\varepsilon>0$, we introduce the associated time subdivision $T_\kappa =N\left[t^\star(\varepsilon) + \kappa \varepsilon\right]$.
    For convenience, we denote $W_{\kappa,\varepsilon}=\macrostack\left(\kappa\varepsilon +\frac{\alpha }{\lambda_ 2 (p-q)}\right)$ and $R_{\kappa,\varepsilon} = \macrostack(t^\varepsilon(\varepsilon)+ \kappa \varepsilon)$. \\

    We also prove ($H_{\kappa,2}$) by induction on $\kappa$. Let us sketch the proof. The hydrodynamic limit on the segment implies that on a $o(N)$ window near $N$, there are $NR_{\kappa,\varepsilon}$.  The strategy consists of slowing down the $NW_{\kappa,\varepsilon}$ first particles by comparing them to a \ZRP with a constant rate. Consequently, we can use the analogous result proven by Labbé \& Lacoin \cite[Proposition 9]{Labb__2019} to control the flux through $N$.   \\
    
     For $\kappa=1$, the result is trivial. Assume ($H_{\kappa,2}$) holds and that {$t^\star(\varepsilon)+ \kappa \varepsilon \leq T_{eq,p,\alpha}$}. Given $\gamma>0$, we introduce three parameters $\gamma_1, \gamma_2, \beta >0$ that we will choose later. \\
    \begin{figure}[h!]
    \label{fig couplage right}
    \centering
    \includegraphics{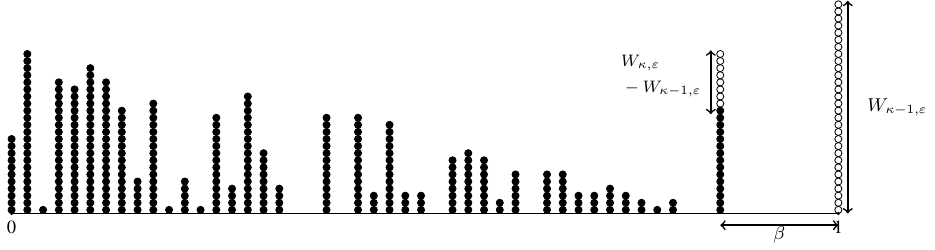}
\caption{The configuration $\tilde \eta$ at $T_\kappa$ zoomed out by a factor $N$.}
\end{figure}
    We define a coloured process $\tilde \eta=\tilde \eta ^b +\tilde \eta^w$ (figure \ref{fig couplage right}) where we denoted the black particles by $\tilde \eta ^b$ and $\tilde \eta ^w$ the white particles. We initialize them at time $T_ \kappa$ by:

\[\tilde{\eta}^b\left(x,T_\kappa\right)=
\left \{
\begin{array}{lr}
    \wedge_{T_\kappa}(x)& \text{if} \quad x \leq x \leq N- \lceil N \beta \rceil,\\
    k - h _{T_ \kappa}(N - \lceil N \beta\rceil -1)- NW_ {\kappa,\varepsilon} & \text{if} \quad x =N- \lceil N \beta \rceil,\\
   0 & \text{otherwise}.
\end{array}
\right.
\]
and 
\[\tilde{\eta}^w\left(x,T_\kappa\right)=
\left \{
\begin{array}{lr}
    0 & \text{if} \quad x \leq x \leq N-\lceil N \beta \rceil,\\
    N[W_ {\kappa,\varepsilon}- W_ {\kappa-1,\varepsilon} + \gamma_1]& \text{if} \quad x = N-\lceil N \beta \rceil,\\
    N (W_{\kappa-1,\varepsilon}- \gamma_1) & \text{if} \quad x= N,\\
   0 & \text{otherwise} .\\
\end{array}
\right.
\]
In particular, on the event 
$$\mathcal{E}=\{S_ {N,k}(T_ \kappa)   \geq N\left[W_{\kappa-1,\varepsilon}-\gamma_1\right]\} \cap \{k - h_{T_ \kappa}\left(N - \lceil \beta N \rceil \rceil\right) \geq N(R_{\kappa,\varepsilon}- \gamma_2)\},$$
the height functions are ordered $$\tilde h_{T_\kappa} ^b +\tilde h_{T_ \kappa }^w \geq h_{T_\kappa},$$
where $\tilde h ^b, \tilde h^w , h_{T_\kappa}$ are the corresponding height functions.\\

Regarding the dynamics, the white particles evolve only on $\llbracket N- \lceil N\beta\rceil,N \rrbracket $ with a modified rate function $\tilde{g}=(pg(1)+q \Bar{g}) \ind_{\cdot\geq 1}$ and a bias $\tilde{p}= \frac{pg(1)}{pg(1)+q \bar g}$. The black particles evolve only on $I_{N-\lceil N\beta\rceil} $ according to a \ZRP$(g,p)$. Let us emphasise that the black and white particles coexist only on the site $N-\lceil N\beta\rceil$. For this particular site, the jump rate for $\tilde \eta ^b$ to the left is $q g(\tilde \eta(N-\lceil N\beta\rceil))$. In other words, only black particles can see the white ones. We claim the following:
\begin{Lemma}
    \label{Lemma couplage pour le stack}
    There exists a coupling which ensures almost surely on $\mathcal{E}$:
    $$
     \forall \tau \geq t \geq T_\kappa, \qquad
    \tilde h_ t ^b + \tilde h_ t ^w \geq  h_t,
    $$
    where $$\tau= \inf\left \{s \geq T_ \kappa, \quad \tilde \eta_s^b \left ( N - \lceil N \beta \rceil  \right )\leq 1\right\}.$$
\end{Lemma}
The exact form of the generator of $(\tilde \eta^b,\tilde \eta ^w) $ and the proof of this lemma are delayed to the appendix.\\
From now on, the proof is similar to its counterpart (Lemma \ref{Lemma convergence of L partial}). We obtain the following inequality:
\begin{equation}
    \label{eq: control S subdivision evenement}
    \begin{aligned}
    &\P{S_{N,k}\left(T_{\kappa +1}\right) \leq N\left[ W_{\kappa,\varepsilon}-\gamma\right]} \\ & \qquad
    \leq  \P{ \left\{\tilde \eta^w_{T_{\kappa +1}}\left(N\right)\leq N\left[ W_{\kappa,\varepsilon}-\gamma\right]\right\}\cap \mathcal{E}} +\P{{\mathcal{E}^c}} +\P{\{\tau \leq N\varepsilon\} \cap \mathcal{E}}.
    \end{aligned}
\end{equation}
The hydrodynamic limit of $(\wedge_t)$ ensures that
$$\lim_{N \rightarrow +\infty} \P{k - h_{T_ \kappa}\left(N - \lceil \beta N \rceil \rceil\right) \geq N(R_{\kappa,\varepsilon}- \gamma_2)\},}=1.$$
Thus, thanks to the induction assumption $H_{\kappa,2}$, the event $\mathcal{E}$ occurs with probability converging to one.\\
If $\tau \leq N\varepsilon$ on $\mathcal{E}$, there was at least $NR_{\kappa, \varepsilon}-N \gamma_2-N W_{\kappa,\varepsilon}-1$ jump to the left from the site $N-\lceil N\beta\rceil$. Thus
$$\P{\{\tau \leq N\varepsilon \} \cap \mathcal{E}} \leq \P{ \mathcal{P}(N \varepsilon \bar g) \geq NR_{\kappa, \varepsilon}-N \gamma_2-N W_{\kappa,\varepsilon}-1 }.$$
By using the convexity of $\macrostack$,
$$
    R_{\kappa,\varepsilon}- W_{\kappa,\varepsilon} \geq R_{0,\varepsilon}-W_{0,\varepsilon} \geq 6 \varepsilon  \bar g q.
    $$
If we fix $\gamma_2 \leq \varepsilon \bar g q$, $\P{\{\tau \leq N\varepsilon\} \cap \mathcal{E}}$ converges to zero.\\

Regarding the third term of (\ref{eq: control S subdivision evenement}), we remark that $\tilde \eta ^w$ behaves like a \ZRP$(\tilde{g},\tilde p)$ on a segment of size $\beta N$ and with $N W_{\kappa,\varepsilon}$ particles. By mapping the process into an \ASEP (Lemma \ref{Lemma: link ZRP-exclusion}),  \cite[Proposition 13]{Labb__2019} implies:
$$
\forall t \geq \frac{(\sqrt{\beta}+\sqrt{W_{\kappa,\varepsilon}-W_{\kappa-1,\varepsilon}+\gamma_1})^2}{pg(1)-q \bar g}, \qquad  \lim_{\substack{N \rightarrow \infty}} \P{ \tilde \eta^w \left(N t +T_ \kappa, N\right) \geq N\left[ W_{\kappa,\varepsilon}-\gamma\right]} =1.
$$
 As $\macrostack$ is convex, we have $$ W_{\kappa,\varepsilon} -W_{\kappa-1,\varepsilon} \leq \varepsilon (\macrostack)^\prime(T_{eq,p,\alpha}) =\varepsilon (p-q)\Phi\left((-\Phi^\prime)^{-1}\left(\frac{-1}{T_{eq,1,\alpha}}\right)\right)$$
 and we assumed $${pg(1)-q \bar g}>({p-q})  \Phi\left((-\Phi^\prime)^{-1}\left(\frac{-1}{T_{eq,1,\alpha}}\right)\right).$$
 We can choose $\beta,\gamma_1$ small enough to verify 
\begin{equation}
    \label{eq: critere comparison preuve right most}
    \varepsilon \geq \frac{(\sqrt{\beta}+\sqrt{W_{\kappa,\varepsilon}-W_{\kappa-1,\varepsilon}+\gamma_1})^2}{pg(1)-q \bar g},
\end{equation}
which concludes the proof up to Lemma \ref{Lemma couplage pour le stack}.
\end{proof}
\begin{remark}
    \label{rem cond explication 2}
    We observe that the flux of \ZRP$(g,p)$ when density blows up converges to $(p-q) \bar g$.  By slowing down the white particles, the maximal flux reachable is the maximal flux of a \ZRP$(\tilde g, \tilde p)$, that is $pg(1)-q \bar g$. So, (\ref{cond sur p et alpha}) means that we impose to the highest flux reached at $(\macrostack)^\prime(T_{eq,p,\alpha})$ to be below the critical level  $pg(1)-q \bar g$.
\end{remark}
\begin{appendix}
\section{Proof of Lemma \ref{Lemma couplage pour le stack}}
Let us first make explicit the generator of the process $(\tilde \eta^w,\tilde \eta^b)$:
 \begin{align*}
     \mathcal{A}^cJ(\tilde \eta^w,\tilde \eta^b)&= \sum_{i=1}^{N-\lceil N \beta \rceil -1 } pg(\eta^b_i) \left[ J(\eta ^w,(\eta^b)^{i,i+1})-J(\eta ^w,\eta^b)\right]\\
     &+ \sum_{i=N-\lceil N \beta \rceil}^{N-1 } pg(1) \ind_{\eta^w_{i}>1} \left[ J((\eta^w)^{i,i+1},\eta^b)-J(\eta ^w,\eta^b)\right]\\
    &+\sum_{i=2}^{N-\lceil N \beta \rceil -1 } qg(\eta^b_i) \left[ J(\eta ^w,(\eta^b)^{i,i-1})-J(\eta ^w,\eta^b)\right] \\
    &+ qg(\eta^w_{N-\lceil N \beta \rceil}+\eta^b_{N-\lceil N \beta \rceil}) \ind_{\eta^b_{N-\lceil N \beta \rceil \geq 1}} \left[ J(\eta ^w,(\eta^b)^{{N-\lceil N \beta \rceil},{N-\lceil N \beta \rceil}-1})-J(\eta ^w,\eta^b)\right]\\
    &+\sum_{i=N-\lceil N \beta \rceil+1}^{N } q\bar g \ind_{\eta^w_{i}>1} \left[ J((\eta^w)^{i,i-1},\eta^b)-J(\eta ^w,\eta^b)\right].
 \end{align*}
\begin{proof}[Proof of Lemma \ref{Lemma couplage pour le stack}]
In order to justify the existence of the coupling, we compare the rate jumps that may inverse the order of the height functions.
\begin{itemize}
    \item The departure site $i < N- \lceil N\beta\rceil $, only black particle can jump.
\begin{itemize}
    \item A jump to the right of a black particle when $ \tilde h ^b (i)= h(i) $. In such case, $\tilde \eta_i^b \leq \eta_ i$ which implies that $g(\tilde \eta_i^b) \leq g(\eta_ i)$.
    \item A jump to the left of $(\wedge_t)$-particle when $ \tilde h ^b (i-1)= h(i-1) $. In such case, $\tilde \eta_i^b \geq \eta_ i$ which implies that $g(\tilde \eta_{i+1}^b) \geq g(\eta_{i+1})$.
\end{itemize}
\item When $i = N- \lceil N\beta\rceil $,
\begin{itemize}
    \item A jump to the right of a white particle when $ \tilde h ^b (i)+\tilde h^w(i)= h(i) $. In such case, $\tilde \eta_i^w \leq \tilde \eta_i^w + \tilde \eta_i^b \leq \eta_ i$ which implies that $g(1)\ind_{\tilde \eta^w_i >1}  \leq g(\eta_ i)$.
    \item A jump to the left of $(\wedge_t)$-particle when $ \tilde h ^b (i-1)= h(i-1) $. In such case, $\tilde \eta_i^w+\tilde \eta_i^b \geq \eta_ i$ which implies that $g(\tilde \eta_i^w+\tilde \eta_i^b )  \geq g(\eta_{i})$. Moreover, the stopping time ensures the presence of a black particle on this site. Thus, this jump cannot occur.
\end{itemize}
\item The departure site $i > N- \lceil N\beta\rceil $, only black particle can jump.
\begin{itemize}
    \item A jump to the right of a black particle when $ \tilde h ^b (i)+\tilde h^w(i)= h(i) $. In such case, $\tilde \eta_i^w \leq \eta_ i$ which implies that $g(1)\ind_{\tilde \eta^w_i >1}  \leq g(\eta_ i)$.
    \item A jump to the left of $(\wedge_t)$-particle when $ \tilde h ^b (i-1)+\tilde h_ t ^w (i-1)= h(i-1) $. In such case, $\tilde \eta_i^w \geq \eta_ i$ which implies that $\bar g  \ind_{\tilde \eta^w_i >1} \geq g(\eta_{i+1})$.
\end{itemize}
\end{itemize}
From this discussion, we conclude the existence of a coupling analogous to \ref{eq: coupling attractivity} and we omit its tedious expression.
\end{proof}
\section{A technical lemma}
For the sake of completeness, we provide proof of the concentration of Poisson variables that we used on multiple occasions.
\begin{Lemma}
    \label{lem concentration poisson}
    Given any $\alpha>0$ and $C>0$, let $(X_i)$ a sequence of independent variable distributed according to $\mathcal{P}(\alpha N)$, then 
    \[
      \lim_{\substack{N \rightarrow \infty}} \P{\max_{i\leq CN}X_i \geq 3\alpha N} =0.
    \]
\end{Lemma}
\begin{proof}
     Let $\gamma= \ln\left[{\frac{\ln{N}}{N\alpha}+1}\right]$, Markov's inequality implies
    \begin{align*}
        \P{\max_{i\leq CN}X_i \geq 3\alpha N} & \leq \frac{\E{e^{\gamma \max_{i\leq CN}X_i}}}{e^{3\gamma\alpha N}}\\
        &\leq C N \frac{\E{e^{\gamma X_1}}}{e^{3\gamma\alpha N}} \leq C N \frac{e^{\alpha N (e^\gamma -1)}}{e^{3\gamma\alpha N}}\\
        & \leq {C}{N^2} e^{-3\alpha N \ln\left[{\frac{\ln{N}}{N\alpha}+1}\right]}.
    \end{align*}
    which converges to zero.
\end{proof}
\end{appendix}
\section*{List of symbols}
\begin{description}
    \item[$N$]{The size of the segment.}
    \item[$k$]{The number of paricles.}
    \item[$g$]{The rate function of the ZRP}
    \item[$\mathcal{A}_N$]{Generator of the \ZRP on the segment $\llbracket 1,N\rrbracket$.}
    \item[$\Omega_{N,k}^0$]{The space of configurations of $k$ particles on the segment $\llbracket 1, N\rrbracket$.}
    \item[$\Omega_{N,k}$]{The space of height functions of configurations in $\Omega_{N,k}^0$.}
    \item[$\pi_N^k$] {The equilibrium measure of the \ZRP on the segment with $k$ particles.}
    \item[$\alpha$]{An asymptotic density.}
    \item[$\Phi$]{The asymptotic flux function.}
    \item[$\Psi$]{The convex conjugate of $\Phi$ if it is convex and of $-\Phi$ if $\Phi$ is concave.}
    \item[$d_{N,k}$]{The worst--case total--variation distance of \ZRP to the equilibrium.}
    \item[$T_{\operatorname{mix}}^ {N,k}$]{Mixing time of the ZRP on the segment with $k$ particles on the segment $\llbracket 1,N\rrbracket$.}
    \item[$U^\alpha$]{The Barron--Jensen solution to (\ref{eq: HJ on R}).}
    \item[$T_{eq,p,\alpha}$]{The macroscopic equilibrium time of $U^\alpha$.}
    \item[$\wedge_{N,k}$]{The configuration where all particles are localized at position $1$.}
    \item[$\vee_{N,k}$]{The configuration where all particles are localized at position $N$.}
    \item[$\ell(\eta)$]{The position of the left--most particle in a configuration $\eta$. }
    \item[$L_{N,k}$]{The position of the left-most particle for \ZRP started from $\wedge_{N,k}$.}
    \item[$S_{N,k}$]{The  mass accumulated at the position $N$ for \ZRP started from $\wedge_{N,k}$.}
    \item[$ \mathscr{L}_p^\alpha$]{The left-most position where a macroscopic mass is observed starting from $\alpha \delta_0$.}
    \item[$ \mathscr{S}_p^\alpha$]{The macroscopic mass accumulated at the position $1$ starting from $\alpha \delta_0$.}
    \item[$\mathcal{P}(\alpha)$]{A Poisson distribution of parameter $\alpha$.}
    \item[$\operatorname{Cont}$]{The set of continuity point of a given function.}
    \item[$f^{\star}$]{The upper semi--continuous envelope of a given function $f$.}
\end{description}

\bibliographystyle{plain}
\bibliography{biblio}
\end{document}